# GENERAL FREQUENTIST PROPERTIES OF THE POSTERIOR PROFILE DISTRIBUTION[1]

By Guang Cheng and Michael R. Kosorok

*Duke University and University of North Carolina at Chapel Hill*

In this paper, inference for the parametric component of a semiparametric model based on sampling from the posterior profile distribution is thoroughly investigated from the frequentist viewpoint. The higher-order validity of the profile sampler obtained in Cheng and Kosorok [*Ann. Statist.* **36** (2008)] is extended to semiparametric models in which the infinite dimensional nuisance parameter may not have a root-$n$ convergence rate. This is a nontrivial extension because it requires a delicate analysis of the entropy of the semiparametric models involved. We find that the accuracy of inferences based on the profile sampler improves as the convergence rate of the nuisance parameter increases. Simulation studies are used to verify this theoretical result. We also establish that an exact frequentist confidence interval obtained by inverting the profile log-likelihood ratio can be estimated with higher-order accuracy by the credible set of the same type obtained from the posterior profile distribution. Our theory is verified for several specific examples.

**1. Introduction.** Semiparametric models have the form $\mathcal{P} = \{P_{\theta,\eta} : (\theta, \eta) \in \Theta \times \mathcal{H}\}$, where $\Theta \subset \mathbb{R}^d$ and $\mathcal{H}$ is an arbitrary subset that is typically infinite dimensional. In this paper, interest will focus on the parametric component $\theta$, while the nonparametric component $\eta$ will be considered a "nuisance parameter." Inference for $\theta$ will be based on semiparametric maximum likelihood estimation via the profile likelihood $pl_n(\theta) = \sup_{\eta \in \mathcal{H}} lik_n(\theta, \eta)$, where $lik_n(\theta, \eta)$ is the full likelihood given $n$ observations. The maximum likelihood estimator for $(\theta, \eta)$ can be expressed as $(\hat{\theta}_n, \hat{\eta}_n)$, where $\hat{\eta}_n = \hat{\eta}_{\hat{\theta}_n}$ and $\hat{\eta}_\theta = \arg\max_{\eta \in \mathcal{H}} lik_n(\theta, \eta)$. We will assume throughout this paper that evaluation of $pl_n(\theta)$ is computationally feasible because of the availability of

Received July 2007; revised August 2007.

[1]Supported in part by Grant CA075142.

procedures such as the stationary point algorithm (as used in [12], e.g.) or the iterative convex minorant algorithm introduced in [7], to find $\hat{\eta}_\theta$ when $\eta$ is a monotone function.

Many of the advantages of using the profile sampler proposed in [14] for inference on $\theta$ are discussed in [14]. The main argument is that direct maximization of the full likelihood and direct computation of the efficient Fisher information function, which often requires tedious evaluation of infinite-dimensional operators that may not have a closed form, can both be avoided completely by using the profile sampler. This follows because the profile sampler yields a first-order correct approximation to the maximum likelihood estimator $\hat{\theta}_n$ and consistent estimation of the efficient Fisher information for $\theta$, even when the nuisance parameter is not estimable at the $\sqrt{n}$ rate.

Another approach to obtaining inference on $\theta$ is through the fully Bayesian procedure, which assigns a prior on both the parameter of interest and the functional nuisance parameter. The first-order valid results in [21] indicate that the marginal semiparametric posterior is asymptotically normal and centered at the corresponding maximum likelihood estimator or posterior mean, with covariance matrix equal to the inverse of the efficient Fisher information. Assigning a prior on $\eta$ can be quite challenging since for some models there is no direct extension of the concept of a Lebesgue dominating measure for the infinite-dimensional parameter set involved [13]. Comparing to the profile sampler procedure, this marginal approach does not circumvent the need to specify a prior on $\eta$, with all of the difficulties that entails. However, we can essentially generate the profile sampler from the marginal posterior of $\theta$ with respect to a certain joint prior on $\psi = (\theta, \eta)$, which is possibly data dependent. For example, we can use a gamma process prior on $\eta$ with jumps at observed event times but not involving $\theta$ in the Cox model with right censored data, see Remark 7 in [3].

The first-order validity of the profile sampler procedure established by [14] is extended to second-order validity in [4] when the infinite-dimensional nuisance parameter achieves the parametric rate. Specifically, higher-order estimates of the maximum profile likelihood estimator and of the efficient Fisher information are obtained in [4]. Moreover, [4] also proves that an exact frequentist confidence interval for the parametric component at level $\alpha$ can be estimated by the $\alpha$-level credible set from the profile sampler with an error of order $O_P(n^{-1})$. Three rather different semiparametric models, the Cox model with right-censored data, the proportional odds model with right-censored data and case-control studies with a missing covariate, are studied in [4]. Such higher-order frequentist accuracy had not previously been established in semiparametric models for any other inferential approach, including the bootstrap. We note that this idea of higher-order accuracy is quite distinct from the concept of second-order efficiency in semiparametric models (see [6, 8]) which we do not consider further in this paper.

A natural question is whether the second-order extension in [4] can be further extended to settings where the nuisance parameter has arbitrary convergence rates, in particular, rates that are slower than the parametric rate. This extension is the key purpose of this paper. Additionally, we generalize the results to allow for multivariate parametric components (only univariate components were permitted in [4]) and also show that another type of confidence interval for $\theta$, obtained by inverting the profile log-likelihood ratio, can also be estimated with higher-order accuracy by the profile sampler.

In this paper, the convergence rate for the nuisance parameter $\eta$ is defined as the largest $r$ that satisfies $\|\hat{\eta}_{\tilde{\theta}_n} - \eta_0\| = O_P(\|\tilde{\theta}_n - \theta_0\| + n^{-r})$, where $\eta_0$ is the true value of $\eta$ and $\|\cdot\|$ is a norm with definition depending on context, that is, for a Euclidean vector $u$, $\|u\|$ is the Euclidean norm, and for an element of the nuisance parameter space $\eta \in \mathcal{H}$, $\|\eta\|$ is some chosen norm on $\mathcal{H}$. In regular semiparametric models, which we can define without loss of generality to be models where the entropy integral converges, $r$ is always larger than $1/4$. Obviously $\hat{\eta}_{\tilde{\theta}_n} \xrightarrow{p} \eta_0$ for any $\tilde{\theta}_n \xrightarrow{p} \theta_0$. We say the nuisance parameter has parametric rate if $r = 1/2$. For instance, the nuisance parameters of the three examples in [4] achieve the parametric rate. More specifically, the nuisance parameter in the Cox model, which is the cumulative hazard function, has the parametric rate under right censored data. However, the convergence rate for the cumulative hazard becomes slower, that is, $r = 1/3$, under current status data. The result is not surprising since current status data cannot provide as much information as right-censored data.

Obviously our results for $r = 1/2$ coincide with the results in [4]. It is also no surprise that the accuracy of the profile sampler is dependent on the convergence rate of the nuisance parameter. The precise error rate for many of the quantities we study is $O_P(M_n(r))$, where we define $M_n(r) = n^{-1/2} + n^{-2r+1/2}$ with support $r > 1/4$. Note that $M_n(r)$ increases in $r$ over the interval $1/4 < r < 1/2$ and is constant for $r \geq 1/2$. Although we cannot yet prove it, we conjecture that this error rate is sharp, in the sense that when the error is multiplied by $M_n^{-1}(r)$, it converges to a nondegenerate random quantity as $n \to \infty$.

Perhaps the most important new result in this paper involves a comparison between an exact, frequentist confidence interval and a credible set for $\theta$ generated from the profile sampler. Specifically, we show that any rectangular credible set for $\theta$ of level $1 - \alpha$ based on the profile sampler is within $O_P(n^{-1/2}M_n(r))$ of an exact, frequentist, rectangular confidence region with coverage $1 - \alpha$. Note that the choice of a one-sided credible set at a given level is not unique when the parameter dimension is $\geq 2$. We also establish higher-order accuracy for the confidence interval obtained by inverting the profile log-likelihood ratio, defined as $PLR_f(\theta) = 2(\log pl_n(\hat{\theta}_n) - \log pl_n(\theta))$.

The next section, Section 2, provides some necessary background material on semiparametric models, least favorable submodels, and empirical processes. The main concepts are illustrated with three examples which will be used throughout the paper. The primary assumptions required for the results of the paper are also presented, along with a key tool for obtaining rates of convergence. In Section 3, second-order asymptotic expansions of the log-profile likelihood are presented. In Section 4, we present the main result of the paper that the confidence interval for the parametric component of a semiparametric model can be approximated by the credible set based on the profile sampler with error of order $O_P(M_n(r))$. In Section 5, we establish that the required assumptions are satisfied for the three previously introduced examples and present some simulation results. Section 6 contains a brief discussion of future research directions, and proofs are given in the Appendix.

**2. Background and assumptions.** We assume the data $X_1, \ldots, X_n$ are i.i.d. throughout the paper. In what follows, we first briefly review the concept of a least favorable submodel. We then present three different examples for which we discuss the forms of the least favorable submodel and related model specifications. Next, we present the model assumptions needed for the remainder of the paper, and, finally, we give a key tool for the rate of convergence calculations needed in later sections.

2.1. *The least favorable submodel.* The *score function* for $\theta$, $\dot{\ell}_{\theta,\eta}$, is defined as the partial derivative w.r.t. $\theta$ of the log-likelihood given $\eta$ is fixed for a single observation. A score function for $\eta_0$ is of the form

$$\frac{\partial}{\partial t}\Big|_{t=0} \log p_{\theta_0,\eta_t}(x) \equiv A_{\theta_0,\eta_0} h(x),$$

where $h$ is a "direction" by which $\eta_t \in \mathcal{H}$ approaches $\eta_0$, running through some index set $H$. $A_{\theta,\eta} : H \mapsto L_2^0(P_{\theta,\eta})$ is the score operator for $\eta$. The *efficient score function* for $\theta$ is defined as $\tilde{\ell}_{\theta,\eta} = \dot{\ell}_{\theta,\eta} - \Pi_{\theta,\eta}\dot{\ell}_{\theta,\eta}$, where $\Pi_{\theta,\eta}\dot{\ell}_{\theta,\eta}$ minimizes the squared distance $P_{\theta,\eta}(\dot{\ell}_{\theta,\eta} - k)^2$ over all functions $k$ in the closed linear space of the score functions for $\eta$ (the "nuisance scores"). A submodel $t \mapsto p_{t,\eta}$ is defined to be *least favorable* at $(\theta, \eta)$ if $\tilde{\ell}_{\theta,\eta} = \partial/\partial t \log p_{t,\eta}$, given $t = \theta$. The inverse of the variance of $\tilde{\ell}_{\theta,\eta}$ is the Crámer–Rao bound for estimating $\theta$ in the presence of the infinite-dimensional nuisance parameter $\eta$, the efficient information matrix $\tilde{I}_{\theta,\eta}$. We also abbreviate $\tilde{\ell}_{\theta_0,\eta_0}$ and $\tilde{I}_{\theta_0,\eta_0}$ with $\tilde{\ell}_0$ and $\tilde{I}_0$, respectively. The direction $h$ along which $\eta_t$ approaches $\eta$ in the least favorable submodel is called "*the least favorable direction.*" An insightful review of least favorable submodels and efficient score functions can be found in Chapter 3 of [11].

The least favorable submodel in this paper will be constructed in the following manner: We first assume the existence of a smooth map from the neighborhood of $\theta$ into the parameter set for $\eta$, of the form $t \mapsto \eta_t(\theta, \eta)$, such that the map $t \mapsto \ell(t, \theta, \eta)(x)$ can be defined as follows:

$$\ell(t, \theta, \eta)(x) = \log lik(t, \eta_t(\theta, \eta))(x), \tag{1}$$

where $t$ and $\theta$ are allowed to be multi-dimensional in this paper, although they both must have the same dimension, and where we require $\eta_\theta(\theta, \eta) = \eta$ for all $(\theta, \eta) \in \Theta \times \mathcal{H}$. We will now illustrate the form of this map for several examples, and the remaining requirements for the map will be presented when the model assumptions are listed later on in this section.

2.2. *Examples.* Three examples with different convergence rates are presented in this subsection. The Cox model with right-censored data, which has a parametric convergence rate, has previously been studied in [4]. Nevertheless, it will be useful to review this example briefly here, although most of the details are given in [4]. The second example, the Cox model with current status data, has a cube-root convergence rate for the nuisance parameter. The last example is the partly linear regression model with normal residual error, where the convergence rate of the nuisance parameter is $n^{-2/5}$ under current status data.

2.2.1. *Example 1. The Cox model with right-censored data.* In the Cox proportional hazards model, the hazard function of the survival time $T$ of a subject with covariate $Z$ is expressed as

$$\lambda(t|z) \equiv \lim_{\Delta \to 0} \frac{1}{\Delta} Pr(t \le T < t + \Delta | T \ge t, Z = z) = \lambda(t) \exp(\theta z), \tag{2}$$

where $\lambda$ is an unspecified baseline hazard function and $\theta$ is a vector including the regression parameters [5]. For the Cox model applied to right-censored failure time data, we observe $X = (Y, \delta, Z)$, where $Y = T \wedge C$, $\delta = I\{T \le C\}$, and $Z \in \mathbb{Z} \subset \mathbb{R}^d$ is a regression covariate. The cumulative hazard function $\Lambda(y) = \int_0^y \lambda(t) \, dt$ is considered the nuisance parameter. The convergence rate of the estimated nuisance parameter is established in Theorem 3.1 of [17], that is, $\|\hat{\Lambda}_{\tilde{\theta}_n} - \Lambda_0\|_\infty = O_P(n^{-1/2} + \|\tilde{\theta}_n - \theta_0\|)$.

Based on the model assumptions specified in Section 5.1 of [4], we can express the likelihood for $(\theta, \eta)$ in the following form:

$$lik(\theta, \Lambda) = (e^{\theta z} \Lambda\{y\} e^{-e^{\theta z} \Lambda(y)})^\delta (e^{-e^{\theta z} \Lambda(y)})^{1-\delta},$$

by replacing $\lambda(y)$ by the point mass $\Lambda\{y\}$. Hence the score functions for $\theta$ and $\Lambda$ can be easily derived as $\dot{\ell}_{\theta, \Lambda}(x) = \delta z - z e^{\theta z} \Lambda(y)$ and $A_{\theta, \Lambda} h(y, \delta, z) =$ $\delta h(y) - e^{\theta z} \int_{[0,y]} h \, d\Lambda$. Again, by the derivations in Section 5.1 of [4], the least favorable direction at $(\theta, \Lambda)$, denoted $h_{\theta,\Lambda}$, can be shown to be

$$h_{\theta,\Lambda}(y) = \frac{E_{\theta,\Lambda} e^{\theta Z} Z 1\{Y \geq y\}}{E_{\theta,\Lambda} e^{\theta Z} 1\{Y \geq y\}}.$$

If we let $h_0$ denote the least favorable direction at the true parameters, the least favorable submodel $\ell(t, \theta, \Lambda)$ has the form $\ell(t, \theta, \Lambda) = \log lik(t, \Lambda_t(\theta, \Lambda))$, where $t \mapsto \Lambda_t(\theta, \Lambda) = \Lambda + (\theta - t)h_0$. Note that we have tacitly swapped the notation for $\eta$ with $\Lambda$ since $\Lambda$ is more widely used in this context.

### 2.2.2. *Example 2. The Cox model with current status data.*

Current status data arises when each subject is observed at a single examination time, $Y$, to determine if an event has occurred. The event time, $T$, cannot be known exactly. If a vector of covariates, $Z$, is also available, then the observed data are $n$ i.i.d. realizations of $X = (Y, \delta, Z) \in R^+ \times \{0, 1\} \times R$, where $\delta = I\{T \leq Y\}$. The model of the conditional hazard given $Z$ is the same as in the previous example. Throughout the remainder of the discussion, we make the following assumptions: $T$ and $Y$ are independent given $Z$. $Z$ lies in a compact set almost surely and the covariance of $Z - E(Z|Y)$ is positive definite, which guarantees the efficient information $\tilde{I}_0$ to be positive definite. $Y$ possesses a Lebesgue density which is continuous and positive on its support $[\sigma, \tau]$, for which the true nuisance parameter $\Lambda_0$ satisfies $\Lambda_0(\sigma-) > 0$ and $\Lambda_0(\tau) < M < \infty$, and this density is continuously differentiable on $[\sigma, \tau]$ with derivative bounded above and bounded below by zero. Under these assumptions the maximum likelihood estimator of $(\theta, \Lambda)$ exists, $\hat{\theta}_n$ is asymptotically efficient and $\|\hat{\Lambda}_n - \Lambda_0\|_{L_2} = O_p(n^{-1/3})$, where $\|\cdot\|_{L_2}$ is the norm on $L_2([\sigma, \tau])$. Note that the conditions on the density of $Y$ ensure that $\|\Lambda - \Lambda_0\|_{L_2}$ is equivalent to $(\int_\sigma^\tau (\Lambda(y) - \Lambda_0(y))^2 \, dF^Y(y))^{1/2}$, where $F^Y(y)$ is the distribution of the observation time $Y$. Moreover, using entropy methods, [17] extends earlier results of [9], showing that $\|\hat{\Lambda}_{\tilde{\theta}_n} - \Lambda_0\|_{L_2} = O_P(\|\tilde{\theta}_n - \theta_0\| + n^{-1/3})$. It is not difficult to derive the log-likelihood

$$\begin{aligned}
(3) \quad \log lik_n(\theta, \Lambda) = \sum_{i=1}^n &\delta_i \log[1 - \exp(-\Lambda(Y_i) \exp(\theta Z_i))] \\
&- (1 - \delta_i) \exp(\theta Z_i)\Lambda(Y_i).
\end{aligned}$$

The score function takes the form $\dot{\ell}_{\theta,\Lambda}(x) = z\Lambda(y)Q(x; \theta, \Lambda)$, where

$$Q(x; \theta, \Lambda) = e^{\theta z}\left[\delta \frac{\exp(-e^{\theta z}\Lambda(y))}{1 - \exp(-e^{\theta z}\Lambda(y))} - (1 - \delta)\right].$$

Inserting a submodel $t \mapsto \Lambda_t$ such that $h(y) = -\partial/\partial t|_{t=0}\Lambda_t(y)$ exists for every $y$ into the log likelihood and differentiating at $t = 0$, we obtain a score function for $\Lambda$ of the form $A_{\theta,\Lambda} h(x) = h(y) Q(x; \theta, \Lambda)$. The linear span of these functions contains $A_{\theta,\Lambda} h$ for all bounded functions $h$ of bounded variation. The efficient score function for $\theta$ is defined as $\tilde{\ell}_{\theta,\Lambda} = \dot{\ell}_{\theta,\Lambda} - A_{\theta,\Lambda} h_{\theta,\Lambda}$ for the vector of functions $h_{\theta,\Lambda}$ minimizing the distance $P_{\theta,\Lambda} \| \dot{\ell}_{\theta,\Lambda} - A_{\theta,\Lambda} h \|^2$, which is also called the least favorable direction. The solution at the true parameter $(\theta_0, \Lambda_0)$ is $h_0(Y)$ defined as follows:

$$(4) \qquad y \mapsto h_0(y) = \Lambda_0(y) h_{00}(y)$$

$$\equiv \Lambda_0(y) \frac{E_{\theta_0 \Lambda_0}(Z Q^2(X; \theta_0, \Lambda_0) | Y = y)}{E_{\theta_0 \Lambda_0}(Q^2(X; \theta_0, \Lambda_0) | Y = y)}.$$

As the formula shows, the vector of functions $h_0(y)$ is unique a.s., and $h_0(y)$ is a bounded function since $Q(x; \theta_0, \Lambda_0)$ is bounded away from zero and infinity. We shall assume the function $y \mapsto h_0(y)$ given by (4) has a version which is differentiable with a bounded derivative on $[\sigma, \tau]$.

The least favorable submodel can be defined as $\ell(t, \theta, \Lambda) = \log lik(t, \Lambda_t(\theta, \Lambda))$, where $\Lambda_t(\theta, \Lambda) = \Lambda + (\theta - t) \phi(\Lambda)(h_{00} \circ \Lambda_0^{-1}) \circ \Lambda$, and $\phi(\cdot)$ is a specially constructed function that smoothly approximates the identity. The function $\Lambda_t(\theta, \Lambda)$ is essentially $\Lambda$ plus a perturbation in the least favorable direction, $h_0$, but its definition is somewhat complicated in order to ensure that $\Lambda_t(\theta, \Lambda)$ really defines a cumulative hazard function within our parameter space, at least for $t$ that is sufficiently close to $\theta$. The details on the construction of the least favorable submodel can be found on page 23 of [16].

### 2.2.3. *Example* 3. *Partly linear normal model with current status data.*
In this example, a continuous outcome $Y$, conditional on the covariates $(W, Z) \in \mathbb{R}^d \times \mathbb{R}$, is modeled as $Y = \theta^T W + k(Z) + \xi$, where $k$ is an unknown smooth function, and $\xi \sim N(0, 1)$. Note that the choice $N(0, 1)$ is needed for model identifiability. We are interested in the regression parameter $\theta$ and consider $k(\cdot)$ to be an infinite-dimensional nuisance parameter. However, the response $Y$ is not observed directly, but only its current status is observed at a random censoring time $C \in \mathbb{R}$. In other words, we observe $X = (C, \Delta, W, Z)$, where $\Delta = 1_{\{Y \le C\}}$. Additionally $(Y, C)$ is assumed to be independent given $(W, Z)$. Although it is not difficult to generalize to multivariate $\theta$, we restrict our attention to univariate $\theta$ in what follows for ease of exposition.

Under the partly linear model, the log-likelihood for a single observation at $X = x \equiv (c, \delta, w, z)$ can be shown to have the form

$$(5) \qquad \log lik_{\theta, k}(x) = \delta \log \{ \Phi(c - \theta w - k(z)) \}$$

$$+ (1 - \delta) \log \{ 1 - \Phi(c - \theta w - k(z)) \},$$

where $\Phi$ is the standard normal distribution. We further assume that the joint distribution for $(C, W, Z)$ is strictly positive and finite. The covariates $(W, Z)$ are assumed to belong to some compact set $\mathcal{W} \times \mathcal{Z} \subset \mathbb{R}^2$. And the random censoring time $C$ is assumed to have support $[l_c, u_c]$, where $-\infty < l_c < u_c < \infty$. In addition, we assume $E[\mathrm{Var}(W|Z)]$ is strictly positive and $Ek(Z) = 0$.

The regression parameter $\theta$ is assumed to belong to some compact set in $\mathbb{R}^1$, and the functional nuisance parameter $k$ is assumed to belong to $\mathcal{O}_2^M \equiv \{f : J_2(f) + \|f\|_\infty < M\}$ for a known $M < \infty$. The $m$th-order Sobolev norm of a function $f$, $J_m(f)$, is defined as $J_m(f) = [\int_{\mathcal{Z}} (f^{(m)}(z))^2\, dz]^{1/2}$. Here, $m$ is a fixed integer and $f^{(j)}$ is the $j$th derivative of $f(\cdot)$ with respect to $z$. The $m$th-order Sobolev class of functions is the class of functions $f$ supported on some compact set on the real line with $J_m(f) < \infty$. Hence the class $\mathcal{O}_2^M$ is trivially the subset of a second-order Sobolev class of functions, and $k \in \mathcal{O}_2^M$ has known upper bound for both its uniform norm and its Sobolev norm. Note that the asymptotic behavior of penalized log-likelihood estimates in this model have been extensively studied in [15].

We now introduce the least favorable submodel. The score function for $\theta$ is $\dot{\ell}_{\theta,k} \equiv w Q(x; \theta, k)$, where

$$Q(X; \theta, k) = (1 - \Delta) \frac{\phi(q_{\theta,k}(X))}{1 - \Phi(q_{\theta,k}(X))} - \Delta \frac{\phi(q_{\theta,k}(X))}{\Phi(q_{\theta,k}(X))}$$

and $q_{\theta,k}(X) = C - \theta W - k(Z)$. Furthermore, by defining $k_t = k + th$ for $h \in \mathcal{O}_2^M$, we can obtain the score function for $k$ in the direction $h$: $A_{\theta,k} h(x) = h(z) Q(x; \theta, k)$. The least favorable direction $h_{\theta,k}$ minimizes $h \mapsto P_{\theta,k} \|\dot{\ell}_{\theta,k} - A_{\theta,k} h\|^2$. By solving the equation $P_{\theta,k}(\dot{\ell}_{\theta,k} - A_{\theta,k} h) A_{\theta,k} h = 0$, we can obtain the solution at the true parameter values:

$$h_0(z) = \frac{E_0(W Q^2(X; \theta, k)|Z = z)}{E_0(Q^2(X; \theta, k)|Z = z)},$$

where $E_0$ is the expectation relative to the true parameters. Thus the least favorable submodel can be constructed as $\ell(t, \theta, k) = \log lik(t, k_t(\theta, k))$, where $k_t(\theta, k) = k + (\theta - t) h_0$.

Note that the above model would be more flexible if we did not require knowledge of $M$. A sieved estimator could be obtained if we replaced $M$ with a sequence $M_n \to \infty$. The theory we propose in this paper will be applicable in this setting, but, in order to maintain clarity of exposition, we have elected not to pursue this more complicated situation here. Another alternative approach is to use penalization. However, this is beyond the scope of the present paper.

2.3. *Assumptions.* We now present the assumptions that will be used throughout the paper, along with some necessary notation. The dependence on $x \in \mathcal{X}$ of the likelihood and score quantities will be largely suppressed for clarity in this section and hereafter.

For the vector $V$, matrix $M$ and tensor $T$, the notation $V_i$, $M_{i,j}$ and $T_{i,j,k}$ indicate its $i$th, $(i,j)$th and $(i,j,k)$th element, respectively. $M^T$ represents the transpose of the matrix $M$. The derivative of the log-likelihood of the least favorable submodel is with respect to the first argument, $t$. The quantities $\dot{\ell}(t,\theta,\eta)$, $\ddot{\ell}(t,\theta,\eta)$ and $\ell^{(3)}(t,\theta,\eta)$ are separately the first, second and third derivative of $\ell(t,\theta,\eta)$ with respect to $t$. For brevity, we denote $\dot{\ell}_0 = \dot{\ell}(\theta_0,\theta_0,\eta_0)$, $\ddot{\ell}_0 = \ddot{\ell}(\theta_0,\theta_0,\eta_0)$ and $\ell_0^{(3)} = \ell^{(3)}(\theta_0,\theta_0,\eta_0)$, where $\theta_0$, $\eta_0$ are the true values of $\theta$ and $\eta$. Of course, $\ddot{\ell}_0(X)$ can also be written as $\dot{\ell}_0(X)$ based on the construction of the least favorable submodel. The quantity $\ell^{(3)}(t,\theta,\eta)$ is a tensor. We thus define $V^T \otimes P\ell^{(3)}(t,\theta,\eta) \otimes V$ as a $d$-dimensional vector whose $i$th element equals $V^T(\partial^2/\partial t^2)(P\dot{\ell}(t,\theta,\eta))_i V$. Similarly $V^T \otimes P\ell^{(3)}(t,\theta,\eta)$ is a square matrix whose $(i,j)$th element is $V^T(\partial/\partial t)(\partial^2/\partial t_i \partial t_j)\ell(t,\theta,\eta)$. We use $\ell_{t_i,t_j,t_k}(t,\theta,\eta)$ to denote $(\partial^3/\partial t_i \partial t_j \partial t_k)\ell(t,\theta,\eta)$. For the derivatives relative to the other two arguments, $\theta$ and $\eta$, we use the following shortened notation: $\ell_\theta(t,\theta,\eta)$ indicates the first derivative of $\ell(t,\theta,\eta)$ with respect to $\theta$. Similarly, $\ell_{t,\theta}(t,\theta,\eta)$ denotes the derivative of $\dot{\ell}(t,\theta,\eta)$ with respect to $\theta$. Also, $\ell_{t,t}(\theta)$ and $\ell_{t,\theta}(\eta)$ indicate the maps $\theta \mapsto \ddot{\ell}(t,\theta,\eta)$ and $\eta \mapsto \ell_{t,\theta}(t,\theta,\eta)$, respectively. Let the random vector $\varrho_n$ denote $\tilde{I}_0^{1/2}(\theta - \hat{\theta}_n)$, and let $\phi_d(\cdot)$ ($\Phi_d(\cdot)$) represent the density (cumulative distribution) of a $d$-dimensional standard normal random variable ($N_d(0,I)$). The notations $\gtrsim$ and $\lesssim$ mean $\geq$, or $\leq$, up to a universal constant. Define $x \vee y$ ($x \wedge y$) to be the maximum (minimum) value of $x$ and $y$. The symbols $\mathbb{P}_n$ and $\mathbb{G}_n \equiv \sqrt{n}(\mathbb{P}_n - P)$ are used for the empirical distribution and the empirical process of the observations, respectively.

We now make the following assumptions in order to achieve the desired second-order asymptotic expansions of the log-profile likelihood (15). The assumption A2 below guarantees that the least favorable submodel passes through $(\theta,\eta)$:

*Regular assumptions*:

A1. $\theta_0 \in \Theta \subset \mathbb{R}^d$, where $\Theta$ is a compact and $\theta_0$ is an interior point of $\Theta$.
A2. $\eta_\theta(\theta,\eta) = \eta$ for any $(\theta,\eta) \in \Theta \times \mathcal{H}$.
A3. $\tilde{I}_0$ is positive definite.

We next describe the smoothness conditions for the least favorable submodel. Clearly, the assumptions B1 and B2 below are separately the smoothness conditions for the Euclidean parameter $(t,\theta)$ and the infinite-dimensional nuisance parameter $\eta$. In principle, these assumptions directly imply the no-bias conditions:

$$\mathbb{P}_n \dot{\ell}(\theta_0, \tilde{\theta}_n, \hat{\eta}_{\tilde{\theta}_n}) = \mathbb{P}_n \tilde{\ell}_0 + O_P(n^{-1/2} + n^{-r} + \|\tilde{\theta}_n - \hat{\theta}_n\|)^2,$$

$$\mathbb{P}_n \ddot{\ell}(\theta_0, \tilde{\theta}_n, \hat{\eta}_{\tilde{\theta}_n}) = P\ddot{\ell}_0 + O_P(n^{-1/2} + n^{-r} + \|\tilde{\theta}_n - \hat{\theta}_n\|),$$

for $\tilde{\theta}_n \xrightarrow{p} \theta_0$, thus making the profile likelihood behave like a standard parametric likelihood asymptotically.

*Smoothness assumptions*:

B1. The maps

$$(t, \theta, \eta) \mapsto \frac{\partial^{l+m}}{\partial t^l \partial \theta^m} \ell(t, \theta, \eta)$$

have integrable envelope functions in $L_1(P)$ in some neighborhood of $(\theta_0, \theta_0, \eta_0)$, for $(l, m) = (0,0), (1,0), (2,0), (3,0), (1,1), (1,2), (2,1)$.

B2. Assume:

$$(6) \qquad \mathbb{G}_n(\dot{\ell}(\theta_0, \theta_0, \hat{\eta}_{\tilde{\theta}_n}) - \dot{\ell}_0) = O_P(M_n(r) + (n^{1/2-r} \vee 1)\|\tilde{\theta}_n - \theta_0\|),$$

$$(7) \qquad P\ddot{\ell}(\theta_0, \theta_0, \eta) - P\ddot{\ell}(\theta_0, \theta_0, \eta_0) = O(\|\eta - \eta_0\|),$$

$$(8) \qquad P\ell_{t,\theta}(\theta_0, \theta_0, \eta) - P\ell_{t,\theta}(\theta_0, \theta_0, \eta_0) = O(\|\eta - \eta_0\|),$$

$$(9) \qquad P\dot{\ell}(\theta_0, \theta_0, \eta) = O(\|\eta - \eta_0\|^2),$$

for $\tilde{\theta}_n \xrightarrow{p} \theta_0$ and all $\eta$ in some neighborhood of $\eta_0$.

There are three approaches to verifying the smoothness assumption (6), which is essentially a continuity modulus of $|\mathbb{G}_n \dot{\ell}(\theta_0, \theta_0, \eta) - \mathbb{G}_n \dot{\ell}_0|$. If the nuisance parameter has parametric convergence rate, we only need to show that the class of functions

$$\left\{ \frac{\dot{\ell}(\theta_0, \theta_0, \eta) - \dot{\ell}_0}{\|\eta - \eta_0\|} : \text{for } \eta \text{ in some neighborhood of } \eta_0 \right\}$$

belongs to a $P$-Donsker class. Alternatively, if the nuisance parameter has the cubic rate, the continuity modulus of the empirical process turns out to be of the order $O_P(n^{-1/6} + n^{1/6}\|\tilde{\theta}_n - \theta_0\|)$, or equivalently $O_P(n^{-1/6} + \|\hat{\eta}_{\tilde{\theta}_n} - \eta_0\|^{1/2})$. The method used to check this condition depends on the norm of the nuisance parameter and the bracketing entropy number of the class of functions $\mathcal{G} = \{\dot{\ell}(\theta_0, \theta_0, \eta)$ for $\eta$ in some neighborhood of $\eta_0\}$. When $\|\cdot\|$ is the $L_2$ norm or one of its dominating norms, we can make use of Lemma 5.13 in [22]. Another approach is to calculate the order of $E_P^* \|\mathbb{G}_n\|_{\mathcal{F}}$, where $\mathcal{F} \equiv \{(\dot{\ell}(\theta_0, \theta_0, \hat{\eta}_{\tilde{\theta}_n}) - \dot{\ell}_0)/(M_n(r) + (n^{1/2-r} \vee 1)\|\tilde{\theta}_n - \theta_0\|)\}$, by the use of Lemma 3.4.2 in [23]. The last two methods will be respectively employed later on in verifying the assumptions for the second and third main examples.

Boundedness of the Fréchet derivatives of the maps $\eta \mapsto \ddot{\ell}(\theta_0, \theta_0, \eta)$ and $\eta \mapsto \ell_{t,\theta}(\theta_0, \theta_0, \eta)$ is sufficient to ensure validity of conditions (7) and (8). Note that $P\ell_{t,\theta}(\theta_0, \theta_0, \eta_0) = 0$ by the following analysis: Fixing $\eta$ and differentiating $P_{\theta,\eta}\dot{\ell}(\theta, \theta, \eta)$ relative to $\theta$ yields $P_{\theta,\eta}\dot{\ell}_{\theta,\eta}\dot{\ell}(\theta, \theta, \eta)^T + P_{\theta,\eta}\ddot{\ell}(\theta, \theta, \eta) +$ $(\partial/(\partial t))|_{t=\theta} P_{\theta,\eta} \dot{\ell}(\theta, t, \eta) = 0$, since $P_{\theta,\eta} \dot{\ell}(\theta, \theta, \eta) = 0$ for every $(\theta, \eta)$, and since we can choose $(\theta, \eta) = (\theta_0, \eta_0)$. One way to verify (9) is to write

$$P\dot{\ell}(\theta_0, \theta_0, \eta) = P\left[\frac{p_0 - p_{\theta_0, \eta}}{p_0}(\dot{\ell}(\theta_0, \theta_0, \eta) - \dot{\ell}(\theta_0, \theta_0, \eta_0))\right]$$
$$- P\left[\dot{\ell}(\theta_0, \theta_0, \eta_0)\left(\frac{p_{\theta_0, \eta} - p_0}{p_0} - A_0(\eta - \eta_0)\right)\right],$$

where $A_0 = A_{\theta_0, \eta_0}$ and $A_{\theta, \eta}$ is the score operator for $\eta$ at $(\theta, \eta)$, for example, the Fréchet derivative of $\log p_{\theta, \eta}$ relative to $\eta$. Thus, if the $L_2$-norm or one of its dominating norms is applied to $\eta$, it suffices to show, under the given regularity conditions, Fréchet differentiability of $\eta \mapsto \dot{\ell}(\theta_0, \theta_0, \eta)$ plus second-order Fréchet differentiability of $\eta \mapsto lik(\theta_0, \eta)$. Note that (9) is naturally satisfied for the semiparametric models with convex linearity, in which $P\dot{\ell}(\theta_0, \theta_0, \eta)$ is exactly zero.

Finally we assume that the following empirical process conditions hold for $(t, \theta, \eta)$ in some neighborhood of their true values:

*Empirical process assumptions*:

C1. There exists some neighborhood $V$ of $(\theta_0, \theta_0, \eta_0)$ in $\Theta \times \Theta \times \mathcal{H}$ such that the classes of functions $\{(\dot{\ell}(t, \theta, \eta))_{i,j}(x) : (t, \theta, \eta) \in V\}$ and $\{(\ell_{t,\theta}(t, \theta, \eta))_{i,j}(x) : (t, \theta, \eta) \in V\}$ are $P$-Donsker and

$$\{(\ell^{(3)}(t, \theta, \eta))_{i,j,k}(x) : (t, \theta, \eta) \in V\}$$

is $P$-Glivenko–Cantelli, for every $i, j, k = 1, \ldots, d$.

One basic method of showing that a class of functions is $P$-Donsker or $P$-Glivenko–Cantelli involves calculating its (bracketing) entropy number. However this verification can be simplified by building up Glivenko–Cantelli (Donsker) classes from other Glivenko–Cantelli (Donsker) classes by employing preservation techniques in Sections 9.3 and 9.4 of [11]. Also, every $P$-Donsker class $\mathcal{F}$ with integrable envelope function is $P$-Glivenko–Cantelli.

2.4. *Rates of convergence.* The estimation accuracy of the profile sampler method depends mainly on the convergence rate of the estimated nuisance parameter, that is, the value of $r$. We now present two useful results, Theorem 1 and Lemma 1 below, that are useful for determining this rate. These results are Theorem 3.2 and Lemma 3.3 in [17], and the proofs can be found therein. Theorem 1 is an extension from general results on M-estimators to semiparametric M-estimators with nuisance parameters. In Theorem 1, $d_\theta^2(\eta, \eta_0)$ may be thought of as the square of a distance, but it is also true for arbitrary functions $\eta \mapsto d_\theta^2(\eta, \eta_0)$. Let $(\Omega, \mathcal{A}, P)$ be an arbitrary probability space and $T : \Omega \mapsto \bar{\mathbb{R}}$ an arbitrary map. Then we use notations $E^*T$ and $O_P^*(1)$ to represent the outer integral of $T$ w.r.t. $P$ and bounded in outer probability, respectively (see page 6 in [23]).

THEOREM 1. *Assume for any given $\theta \in \Theta_n$, that $\hat{\eta}_\theta$ satisfies $\mathbb{P}_n m_{\theta, \hat{\eta}_\theta} \geq \mathbb{P}_n m_{\theta, \eta_0}$ for given measurable functions $x \mapsto m_{\theta, \eta}(x)$. Assume conditions* (10) *and* (11) *below hold for every $\theta \in \Theta_n$, every $\eta \in \mathcal{V}_n$ and every $\varepsilon > 0$:*

$$
(10) \qquad P(m_{\theta, \eta} - m_{\theta, \eta_0}) \lesssim -d_\theta^2(\eta, \eta_0) + \|\theta - \theta_0\|^2,
$$

$$
(11) \qquad E^* \sup_{\theta \in \Theta_n, \eta \in \mathcal{V}_n, \|\theta - \theta_0\| < \varepsilon, d_\theta(\eta, \eta_0) < \varepsilon} |\mathbb{G}_n(m_{\theta, \eta} - m_{\theta, \eta_0})| \lesssim \phi_n(\varepsilon).
$$

*Suppose that* (11) *is valid for functions $\phi_n$ such that $\delta \mapsto \phi_n(\delta)/\delta^\alpha$ is decreasing for some $\alpha < 2$ and sets $\Theta_n \times \mathcal{V}_n$ such that $P(\tilde{\theta} \in \Theta_n, \hat{\eta}_{\tilde{\theta}} \in \mathcal{V}_n) \to 1$. Then $d_{\tilde{\theta}}(\hat{\eta}_{\tilde{\theta}}, \eta_0) \leq O_P^*(\delta_n + \|\tilde{\theta} - \theta_0\|)$ for any sequence of positive numbers $\delta_n$ such that $\phi_n(\delta_n) \leq \sqrt{n}\delta_n^2$ for every $n$.*

Lemma 1 below is useful for verifying the continuity modulus condition (11) for the empirical process. Define $\mathcal{S}_\delta = \{x \mapsto m_{\theta, \eta}(x) - m_{\theta, \eta_0}(x) : d_\theta(\eta, \eta_0) < \delta, \|\theta - \theta_0\| < \delta\}$ and

$$
(12) \qquad K(\delta, \mathcal{S}_\delta, L_2(P)) = \int_0^\delta \sqrt{1 + H_B(\varepsilon, \mathcal{S}_\delta, L_2(P))} \, d\varepsilon,
$$

where $H_B$ denotes the log of the bracketing entropy number.

LEMMA 1. *Suppose that the functions $(x, \theta, \eta) \mapsto m_{\theta, \eta}(x)$ are uniformly bounded for $(\theta, \eta)$ ranging over a neighborhood of $(\theta_0, \eta_0)$ and that*

$$
(13) \qquad P(m_{\theta, \eta} - m_{\theta_0, \eta_0})^2 \lesssim d_\theta^2(\eta, \eta_0) + \|\theta - \theta_0\|^2.
$$

*Then condition* (11) *is satisfied for any functions $\phi_n$ such that*

$$
\phi_n(\delta) \geq K(\delta, \mathcal{S}_\delta, L_2(P)) \left( 1 + \frac{K(\delta, \mathcal{S}_\delta, L_2(P))}{\delta^2 \sqrt{n}} \right).
$$

*Consequently, we may replace $\phi_n(\delta)$ with $K(\delta, \mathcal{S}_\delta, L_2(P))$ in the conclusion of the previous theorem.*

**3. Second-order asymptotic expansion.** In this section, second-order asymptotic expansions of the log-profile likelihood and maximum likelihood estimator are derived. Their second-order accuracy is proven to be dependent on the order of the convergence rate of the nuisance parameter through the rate function $M_n(r)$ given in the Introduction. Note that the smallest order of $O_P(M_n(r))$, $O_P(n^{-1/2})$, is achieved when the nuisance parameter has parametric or faster rate by the truncation property of the function $M_n(r)$. The assumptions in Section 2 are assumed throughout.

THEOREM 2. *If $\tilde{\theta}_n$ satisfies $(\tilde{\theta}_n - \hat{\theta}_n) = o_P(1)$, then*

$$(14) \qquad \sqrt{n}(\hat{\theta}_n - \theta_0) = \frac{1}{\sqrt{n}} \sum_{i=1}^{n} \tilde{I}_0^{-1} \tilde{\ell}_0(X_i) + O_P(M_n(r)),$$

$$\log pl_n(\tilde{\theta}_n) = \log pl_n(\hat{\theta}_n) - \tfrac{1}{2} n(\tilde{\theta}_n - \hat{\theta}_n)^T \tilde{I}_0(\tilde{\theta}_n - \hat{\theta}_n)$$
$$(15) \qquad\qquad\qquad + O_P(g_r(\|\tilde{\theta}_n - \hat{\theta}_n\|)),$$

*where $g_r(w) \equiv (nw^3 + n^{1-r}w^2 + n^{1-2r}w + n^{-2r+1/2})1\{1/4 < r < 1/2\} + (nw^3 + n^{-1/2})1\{r \geq 1/2\}$.*

REMARK 1. Under regularity conditions, the counterpart of (14) in fully parametric models has error of order $O_P(n^{-1/2})$, which agrees with $O_P(M_n(r))$ when $r \geq 1/2$. Thus, we achieve the parametric bound in semiparametric models only when the nuisance parameter obtains the parametric rate. We also observe a monotonic increase in the error rate as $r$ decreases toward $1/4$.

The asymptotic quadratic expansion (15) can be used to construct an estimator of the standard error of $\hat{\theta}_n$. The estimator is the following "discretized" version of the observed profile information matrix, $\hat{I}_n$, which is the derivative of the profile likelihood (see [17]):

$$(16) \qquad \hat{I}_n(v) \equiv -2 \frac{\log pl_n(\hat{\theta}_n + s_n v) - \log pl_n(\hat{\theta}_n)}{n s_n^2},$$

where direction $v \in \mathbb{R}^d$ and step size $s_n \to 0$. The expansion (15) implies

$$(17) \qquad v^T \tilde{I}_0 v = \hat{I}_n(v) + O_P(h_r(|s_n|)),$$

where $h_r(|s_n|) = g_r(|s_n|)/ns_n^2$. By straightforward analysis, the smallest order of the error term in (17) is $O_P(n^{-r})$ by setting the step size $s_n = O_p(n^{-r})$ and $s_n^{-1} = O_P(n^r)$ when $1/4 < r < 1/2$. However, when $r \geq 1/2$, the smallest order of error in (17) stabilizes at $O_P(n^{-1/2})$ by setting the step size to $s_n = O_p(n^{-1/2})$ and $s_n^{-1} = O_P(n^{1/2})$. In other words, $\hat{I}_n$ can only be a $\sqrt{n}$ consistent estimator of $\tilde{I}_0$ when the convergence rate of the nuisance parameter is faster than or equal to the parametric rate.

The above analysis also leads to good discretized estimators for each element in $\tilde{I}_0$. For instance, with $e_i$ denoting the $i$th unit vector in $\mathbb{R}^d$, we can deduce

$$(\hat{I}_n(e))_{i,j} = -\frac{\log pl_n(\hat{\theta}_n + e_i s_n + e_j s_n) + \log pl_n(\hat{\theta}_n)}{n s_n^2}$$
$$(18) \qquad\qquad + \frac{\log pl_n(\hat{\theta}_n + e_i s_n) + \log pl_n(\hat{\theta}_n + e_j s_n)}{n s_n^2},$$

$$(19) \qquad (\tilde{I}_0)_{i,j} = (\hat{I}_n(e))_{i,j} + O_P(h_r(|s_n|)).$$

**4. Main results and implications.** We now present the main results on the posterior profile distribution. Let $\tilde{P}_{\theta|\tilde{X}}$ be the posterior profile distribution of $\theta$ with respect to the prior $\rho(\theta)$ given the data $\tilde{X} = (X_1, \ldots, X_n)$. Define $\Delta_n(\theta) = n^{-1}\{\log pl_n(\theta) - \log pl_n(\hat{\theta}_n)\}$. We now present the first main result:

THEOREM 3. *Assume the assumptions of Section 2 and also that*

$$(20) \qquad \Delta_n(\tilde{\theta}_n) = o_P(1) \qquad \text{implies that } \tilde{\theta}_n = \theta_0 + o_P(1),$$

*for any sequence $\tilde{\theta}_n \in \Theta$. If proper prior $\rho(\theta_0) > 0$ and $\rho(\cdot)$ has a continuous and finite first-order derivative in some neighborhood of $\theta_0$, then*

$$(21) \qquad \sup_{\xi \in \mathbb{R}^d} |\tilde{P}_{\theta|\tilde{X}}(\sqrt{n}\tilde{I}_0^{1/2}(\theta - \hat{\theta}_n) \leq \xi) - \Phi_d(\xi)| = O_P(M_n(r)).$$

REMARK 2. Based on the conclusions of Theorem 3, we know that the $[1 - \alpha + O_P(M_n(r))]$th one-sided and two-sided credible sets for vector $\theta$ from the profile sampler are $(-\infty, \hat{\theta}_n + n^{-1/2}\tilde{I}^{-1/2}z_{1-\alpha}]$ and $[\hat{\theta}_n - n^{-1/2}\tilde{I}^{-1/2}z_{1-\alpha/2}, \hat{\theta}_n + n^{-1/2}\tilde{I}^{-1/2}z_{1-\alpha/2}]$, respectively, where $z_\alpha$ is a standard normal $\alpha$th quantile for $d$-dimensions and $\tilde{I}$ can be either $\tilde{I}_0$ or $\hat{I}_n$.

The following two corollaries provide several interesting additional second-order properties of the profile sampler:

COROLLARY 1. *Assume the conditions of Theorem 3, and let $f_n(\cdot)$ be the posterior profile density of $\sqrt{n}\varrho_n$ relative to the prior $\rho(\theta)$. Then*

$$(22) \qquad f_n(\xi) = \phi_d(\xi) + O_P(M_n(r)).$$

COROLLARY 2. *Under the conditions of Theorem 3 and recalling that $\varrho_n = \tilde{I}_0^{1/2}(\theta - \hat{\theta}_n)$, we have that if $\theta$ has finite second absolute moment, then*

$$(23) \qquad \hat{\theta}_n = \tilde{E}_{\theta|\tilde{X}}(\theta) + O_P(n^{-1/2}M_n(r)),$$

$$(24) \qquad \tilde{I}_0 = n^{-1}(\tilde{\mathrm{Var}}_{\theta|\tilde{X}}(\theta))^{-1} + O_P(M_n(r)),$$

*where $\tilde{E}_{\theta|\tilde{X}}(\theta)$ and $\tilde{\mathrm{Var}}_{\theta|\tilde{X}}(\theta)$ are the posterior profile mean and posterior profile covariance matrix, respectively.*

REMARK 3. The posterior moments in Corollary 2 are with respect to the posterior profile distribution. Thus we can estimate $\hat{\theta}_n$ with the mean of the profile sampler. Similarly, (each element of) the efficient information matrix can be estimated by (the corresponding element of) the inverse of the covariance matrix of the profile sampler with an error of order $O_P(M_n(r))$. Clearly, a faster convergence rate of the nuisance parameter leads to higher estimation accuracy. We can generalize the arguments used in the proof of Corollary 2 to obtain general results on the posterior moments. For simplicity, assume $\theta$ is one dimensional. Then, provided $\int_{-\infty}^{+\infty} |\theta|^\beta \rho(\theta)\,d\theta < \infty$, we have $\tilde{E}_{\theta|\tilde{X}}\varrho_n^\beta = n^{-\beta/2}EU^\beta + O_P(n^{-(\beta+1)/2} + n^{(-2r+1)-(\beta+1)/2})$, where $\tilde{E}_{\theta|\tilde{X}}\varrho_n^\beta$ is the $\beta$th posterior moment of $\varrho_n$ and $U \sim N(0,1)$.

REMARK 4. We now have two approaches to estimating the efficient information matrix $\tilde{I}_0$. One approach is by numerical analysis as given in (19). Another approach is presented in Corollary 2 as an estimate from the posterior distribution. We prefer estimating $\tilde{I}_0$ with (24) using the profile sampler procedure in semiparametric models with $r \geq 1/2$ since this avoids the issue of choosing the step size in (17) or (24). However for models with $r < 1/2$, the numerical differentiation approach may be worthwhile because of the smaller error rate that may be obtained using (17).

Combining (14) and (23), we obtain

$$\sqrt{n}(\tilde{E}_{\theta|\tilde{X}}(\theta) - \theta_0) = \frac{1}{\sqrt{n}}\sum_{i=1}^n \tilde{I}_0^{-1}\tilde{\ell}_0(X_i) + O_P(M_n(r)).$$

The range of $r$ implies that the mean value of the profile sampler is essentially a semiparametric efficient estimator of $\theta$ even when the nuisance parameter has a slower convergence rate. A similar conclusion appears to hold for other estimators of $\hat{\theta}_n$ based on the profile sampler, including multivariate generalizations of the median.

The second main result is expressed in the following Theorem 4. An $\alpha$th quantile of the posterior profile distribution is any quantity $\tau_{n\alpha} \in \mathbb{R}^d$ that satisfies $\tau_{n\alpha} = \inf\{\xi : \tilde{P}_{\theta|\tilde{X}}(\theta \leq \xi) \geq \alpha\}$, where $\xi$ is an infimum over the given set only if there does not exist a $\xi_1 < \xi$ in $\mathbb{R}^d$ such that $\tilde{P}_{\theta|\tilde{X}}(\theta \leq \xi_1) \geq \alpha$. Because of the assumed smoothness of both the prior and the likelihood in our setting, we can, without loss of generality, assume $\tilde{P}_{\theta|\tilde{X}}(\theta \leq \tau_{n\alpha}) = \alpha$. We can also define $\kappa_{n\alpha} = \sqrt{n}(\tau_{n\alpha} - \hat{\theta}_n)$, that is, $\tilde{P}_{\theta|\tilde{X}}(\sqrt{n}(\theta - \hat{\theta}_n) \leq \kappa_{n\alpha}) = \alpha$. Note that neither $\tau_{n\alpha}$ nor $\kappa_{n\alpha}$ are unique if the dimension for $\theta$ is larger than one. Nevertheless, the following theorem ensures that for each choice of $\kappa_{n\alpha}$ there exists a unique $\hat{\kappa}_{n\alpha}$ based on the data such that $P(\sqrt{n}(\hat{\theta}_n - \theta_0) \leq \hat{\kappa}_{n\alpha}) = \alpha$ and $\|\hat{\kappa}_{n\alpha} - \kappa_{n\alpha}\| = O_P(M_n(r))$:

THEOREM 4. Under the conditions of Theorem 3 and assuming that $\tilde{\ell}_0(X)$ has finite third moment with a nondegenerate distribution, then there exists a $\hat{\kappa}_{n\alpha}$ based on the data such that $P(\sqrt{n}(\hat{\theta}_n - \theta_0) \leq \hat{\kappa}_{n\alpha}) = \alpha$ and $\hat{\kappa}_{n\alpha} - \kappa_{n\alpha} = O_P(M_n(r))$ for each chosen $\kappa_{n\alpha}$.

REMARK 5. Clearly, a faster convergence rate of the nuisance parameter leads to a more accurate estimate of the confidence interval when $1/4 < r < 1/2$. The profile sampler procedure can provide the best estimate for the boundary of the confidence interval in semiparametric models when $r \geq 1/2$. We conjecture that the product of $\sqrt{n} I\{r \geq 1/2\} + n^{2r-1/2} I\{1/4 < r < 1/2\}$ and the $O_P(M_n(r))$ term in Theorem 4 converges to the product of two different nontrivial but uniformly integrable Gaussian processes. Thus we believe the convergence rate in Theorem 4 is optimal.

Theorem 4 states that the Wald-type confidence interval can be approximated by the credible set of the same type based on the profile sampler with error of order $O_P(M_n(r))$. In other words, the boundary of a one-sided confidence interval for $\theta$ at level $\alpha$ can be estimated by the $\alpha$th quantile of the profile sampler with error of order $O_P(n^{-1/2}M_n(r))$. Similar conclusions also hold for the confidence interval obtained by inverting the profile likelihood ratio, as will be shown in Theorem 5 below.

The profile likelihood ratio in the frequentist and Bayesian set-up is separately defined as $PLR_f(\theta_0) = 2(\log pl_n(\hat{\theta}_n) - \log pl_n(\theta_0))$ and $PLR_b(\theta) = 2(\log pl_n(\hat{\theta}_n) - \log pl_n(\theta))$. Thus $\chi_b^{n\alpha}$ is defined by $\chi_b^{n\alpha} = \inf\{\xi : \tilde{P}_{\theta|\tilde{X}}(PLR_b(\theta) \leq \xi) \geq \alpha\}$. As argued previously, we can, without loss of generality, assume that $\tilde{P}_{\theta|\tilde{X}}(PLR_b(\theta) \leq \chi_b^{n\alpha}) = \alpha$. The following theorem ensures that there exists a $\chi_f^{n\alpha}$ based on the data such that $P(PLR_f(\theta_0) \leq \chi_f^{n\alpha}) = \alpha$ and $\chi_f^{n\alpha} - \chi_b^{n\alpha} = O_P(M_n(r))$:

THEOREM 5. Under the conditions of Theorem 4, there exists a $\chi_f^{n\alpha}$ based on the data such that $P(PLR_f(\theta_0) \leq \chi_f^{n\alpha}) = \alpha$ and $\chi_f^{n\alpha} - \chi_b^{n\alpha} = O_P(M_n(r))$.

REMARK 6. The corresponding $\alpha$-level confidence interval and credible set obtained by inverting the profile likelihood ratio can be expressed as $C_f^{n\alpha}(X) = \{\theta \in \Theta : PLR_f(\theta) \leq \chi_f^{n\alpha}\}$ and $C_b^{n\alpha}(X) = \{\theta \in \Theta : PLR_b(\theta) \leq \chi_b^{n\alpha}\}$, respectively. Moreover, the proof of Theorem 5 implies that $\chi_b^{n\alpha} = \chi_{d,\alpha}^2 + O_P(M_n(r))$ and $\chi_b^{n\alpha} = \chi_{d,\alpha}^2 + O_P(M_n(r))$, where $\chi_{d,\alpha}^2$ denotes the $\alpha$th quantile of central chi-square distribution with degree of freedom $d$. Theorem 5 also implies that $\tilde{P}_{\theta|\tilde{X}}(PLR_b(\theta) \leq \chi_{d,\alpha}^2) = \alpha + O_P(M_n(r))$. Thus it appears in this instance that not much is gained by using the posterior profile sampler to calibrate the likelihood ratio confidence interval instead of simply using $\chi_{d,\alpha}^2$.

**5. Examples.** This section illustrates the practicality of the stated conditions by verifying that these assumptions are satisfied for each of the three examples introduced in Section 2. Some simulation results about the Cox regression model are also presented.

5.1. *The Cox model with right-censored data.* Note that this example was considered fully in [4], but we include some of the main ideas here for completeness. We first verify the smoothness conditions B1 and the empirical processes assumptions C1. Under regular conditions, B1 can be easily satisfied since the maps $(t, \theta, \eta) \mapsto (\partial^{l+m}/\partial t^l \theta^m)\ell(t, \theta, \eta)$, whose forms can be found in [4], are uniformly bounded around $(\theta_0, \theta_0, \Lambda_0)$. Notice that the functions $y \mapsto h_0(y)$, $y \mapsto \Lambda_t(y)$ and $z \mapsto \exp(zt)$ for $(t, \theta, \Lambda)$ in the assumed neighborhood of the true values are $P$-Donsker. Thus we can verify C1 by repeatedly employing the Lipschitz continuity preservation property of Donsker classes. The remaining smoothness conditions B2 and condition (20) are separately verified by Lemmas 2 and 3 of [4].

5.2. *The Cox model with current status data.* In this section we verify the regularity conditions for the Cox model with current status data as well as present a small simulation study to gain insight into the moderate sample size agreement with the asymptotic theory.

5.2.1. *Verification of conditions.* We can verify that $\ell(t, \theta, \Lambda)$ defined in Section 2.2.2 above is indeed the least favorable submodel since $\dot\ell(t, \theta, \Lambda) = (z\Lambda_t(\theta, \Lambda)(y) - \phi(\Lambda(y))h_{00} \circ \Lambda_0^{-1} \circ \Lambda(y))Q(x; t, \Lambda_t(\theta, \Lambda))$, evaluated at $t = \theta = \theta_0$ and $\Lambda = \Lambda_0$, is the efficient score function $(z\Lambda_0(y) - h_0(y))Q(x; \theta_0, \Lambda_0)$. Note that we extend the domain of the function $u \mapsto \Lambda_0^{-1}(u)$ to all of $[0, \infty)$ by assigning the value $\sigma$ to all $u \in [0, \Lambda(\sigma))$ and the value $\tau$ to all $u > \Lambda(\tau)$. Substituting $\theta = t$ and $\Lambda = \Lambda_t(\theta, \Lambda)$ in (3) and differentiating with respect to $t$ and $\theta$, we obtain,

$$\dot\ell(t, \theta, \Lambda)(x) = (z\Lambda_t + \dot\Lambda_t)Q(x; t, \Lambda_t),$$

$$\ddot\ell(t, \theta, \Lambda)(x) = \frac{\partial^2 lik(t, \Lambda_t(\theta, \Lambda))/\partial t^2}{lik(t, \Lambda_t(\theta, \Lambda))} - \dot\ell^2(t, \theta, \Lambda),$$

where

$$\frac{\partial^2 lik(t, \Lambda_t(\theta, \Lambda))/\partial t^2}{lik(t, \Lambda_t(\theta, \Lambda))} = Q(x; t, \Lambda_t) \times [z^2\Lambda_t + 2\dot\Lambda_t - e^{tz}(z\Lambda_t + \dot\Lambda_t)^2].$$

Note that $\dot\ell(t, \theta, \Lambda)$ can be written as follows:

$$\dot\ell(t, \theta, \Lambda) = \left[z - \frac{\phi(\Lambda)(y)}{\Lambda_t(\theta, \Lambda)}h_{00} \circ \Lambda_0^{-1} \circ \Lambda(y)\right]\Lambda_t(\theta, \Lambda)(y)Q(x; t, \Lambda_t(\theta, \Lambda)),$$

and the map $u \mapsto ue^{-u}/(1 - e^{-u})$ is bounded and Lipschitz on $[0, \infty)$. Thus we can write $\Lambda(y)Q(x; t, \Lambda) = \psi(e^{tz}, \Lambda(y))$, where the function $\psi$ is bounded and Lipschitz in each argument. Next, note that

$$-\frac{\dot\Lambda_t}{\Lambda_t} = \frac{\phi(\Lambda)h_{00} \circ \Lambda_0^{-1} \circ \Lambda}{\Lambda_t(\theta, \Lambda)} = \frac{(\phi(\Lambda)/\Lambda)h_{00} \circ \Lambda_0^{-1} \circ \Lambda}{1 + (\theta - t)(\phi(\Lambda)/\Lambda)h_{00} \circ \Lambda_0^{-1} \circ \Lambda}.$$

Combining this with the facts that the function $\phi(y)/y$ is bounded and $h_{00} \circ \Lambda_0^{-1}$ is bounded by assumption, we obtain that $\dot{\ell}(t, \theta, \Lambda)$ is bounded. Clearly, $\ddot{\ell}(t, \theta, \Lambda)$ is also uniformly bounded based on the following equation:

$$\frac{\partial^2 lik(t, \Lambda_t(\theta, \Lambda))/\partial t^2}{lik(t, \Lambda_t(\theta, \Lambda))} = Q(x; t, \Lambda_t)\Lambda_t$$
$$\times \left( z^2 + 2\frac{\dot{\Lambda}_t}{\Lambda_t} - e^{tz}\Lambda_t\left( z^2 + 2z\frac{\dot{\Lambda}_t}{\Lambda_t} + \left(\frac{\dot{\Lambda}_t}{\Lambda_t}\right)^2 \right) \right).$$

By similar analysis, $\ell_{t,\theta}(t, \theta, \Lambda)$, $\ell^{(3)}(t, \theta, \Lambda)$, $\ell_{t,t,\theta}(t, \theta, \Lambda)$ and $\ell_{t,\theta,\theta}(t, \theta, \Lambda)$, whose concrete forms can be found in [3], are also uniformly bounded for all $t$ sufficiently close to $\theta$ and all $\Lambda$ varying over the parameter space.

We next verify assumption C1. Recall that $\Lambda(y)Q(x; t, \Lambda) = \psi(e^{tz}, \Lambda(y))$, where the function $\psi$ is bounded and Lipschitz in each argument. Thus, since the classes of functions $z \mapsto e^{tz}$ and $y \mapsto \Lambda(y)$ are Donsker, so is the class of functions $x \mapsto \Lambda(y)Q(x; t, \Lambda)$. Note that

$$\frac{\phi(\Lambda)}{\Lambda_t(\theta, \Lambda)} = \frac{\varsigma(\Lambda)}{1 + (\theta - t)\varsigma(\Lambda)\upsilon(\Lambda)} \equiv \chi(\Lambda),$$

where $\varsigma(\Lambda) = \phi(\Lambda)/\Lambda$ and $\upsilon(\Lambda) = h_{00} \circ \Lambda_0^{-1} \circ \Lambda$, and both $\varsigma(\Lambda)$ and $\upsilon(\Lambda)$ are Lipschitz according to the assumptions. Hence $\chi(\Lambda)$ is also Lipschitz in $\Lambda$. Thus the class of functions $\dot{\ell}(t, \theta, \Lambda)$ with $(t, \theta)$ varying over a small neighborhood of $(\theta_0, \theta_0)$ and $\Lambda$ ranging over all nondecreasing cadlag functions with domain $[\sigma, \tau]$ and range $[0, M]$ can be seen to be a Donsker class. By repeated application of the above techniques to $(\partial^2 lik(t, \Lambda_t(\theta, \Lambda))/\partial t^2)/lik(t, \Lambda_t(\theta, \Lambda))$ we know the class of functions $\ddot{\ell}(t, \theta, \Lambda)$ is also Donsker. Similarly, the classes of functions $\ell_{t,\theta}(t, \theta, \Lambda)$ and $\ell^{(3)}(t, \theta, \Lambda)$ with $(t, \theta)$ varying over a small neighborhood of $(\theta_0, \theta_0)$ and $\Lambda$ ranging over all nondecreasing cadlag functions with domain $[\sigma, \tau]$ and range $[0, M]$ can be shown to be Donsker. Moreover, $\ell^{(3)}(t, \theta, \Lambda)$ is automatically $P$-Glivenko–Cantelli since it is uniformly bounded based on the previous analysis. The following lemmas verify the remaining assumptions:

LEMMA 2.   *Under the above set-up for the Cox model with current status data, assumption* B2 *is satisfied.*

LEMMA 3.   *Under the above set-up for the Cox model with current status data, condition* (20) *is satisfied.*

5.2.2. *Simulation study.*   In this subsection, we conducted simulations in two semiparametric models with different convergence rates, that is, Cox regression with right-censored data and Cox regression with current status data. The contrast of the above two simulations agrees with our theoretical results that the accuracy of inferences based on the profile sampler is higher in semiparametric models with faster convergence rates.

In what follows, the simulations are run for various sample sizes under a Lebesgue prior. For each sample size, 500 datasets were analyzed. The event times were generated from (2) with one covariate $Z \sim U[0,1]$. The regression coefficient is $\theta = 1$ and $\Lambda(t) = \exp(t) - 1$. The censoring time $C \sim U[0, t_n]$, where $t_n$ was chosen such that the average effective sample size over 500 samples is approximately $0.9n$. For each dataset, Markov chains of length 20,000 with a burn-in period of 5,000 were generated using the Metropolis algorithm. The jumping density for the coefficient was normal with current iteration and variance tuned to yield an acceptance rate of 20%–40%. The approximate variance of the estimator of $\theta$ was computed by numerical differentiation with step size proportional to $n^{-1/2}$ ($n^{-1/3}$) for right-censored data (current status data) according to (16).

Table 1 (Table 2) summarizes the results from the simulations of Cox regression with right-censored data (current status data) giving the average across 500 samples of the maximum likelihood estimate (MLE), mean of the profile sampler (CM), estimated standard errors based on MCMC ($\text{SE}_\text{M}$), estimated standard errors based on numerical derivatives ($\text{SE}_\text{N}$) and

TABLE 1
*Cox regression with right-censored data ($\theta_0 = 1$ and 500 samples)*

| $n$ | $n|\text{MLE} - \text{CM}|$ | $\sqrt{n}|\text{SE}_\text{M} - \text{SE}_\text{N}|$ | $n|\text{L}_\text{M} - \text{L}_\text{N}|$ | $n|\text{U}_\text{M} - \text{U}_\text{N}|$ |
|---|---|---|---|---|
| 50 | 0.3062 | 0.2270 | 0.1809 | 1.1212 |
| 100 | 0.2587 | 0.0311 | 0.5987 | 0.1301 |
| 200 | 0.3218 | 0.0279 | 0.4810 | 0.5253 |
| 500 | 0.2017 | 0.2080 | 0.7524 | 0.3518 |

TABLE 2
*Cox regression with current status data ($\theta_0 = 1$ and 500 samples)*

| $n$ | $n^{2/3}|\text{MLE} - \text{CM}|$ | $n^{2/6}|\text{SE}_\text{M} - \text{SE}_\text{N}|$ | $n^{2/3}|\text{L}_\text{M} - \text{L}_\text{N}|$ | $n^{2/3}|\text{U}_\text{M} - \text{U}_\text{N}|$ |
|---|---|---|---|---|
| 50 | 0.4438 | 0.6144 | 3.1799 | 5.6550 |
| 100 | 0.6506 | 0.4071 | 0.6162 | 1.0082 |
| 200 | 0.7729 | 0.3284 | 0.4617 | 0.8071 |
| 500 | 0.7559 | 0.1611 | 0.1071 | 1.3103 |

$n$, sample size; MLE, maximum likelihood estimator; CM, empirical mean; $\text{SE}_\text{M}$, estimated standard errors based on MCMC; $\text{SE}_\text{N}$, estimated standard errors based on numerical derivatives; $\text{L}_\text{M}$ ($\text{U}_\text{M}$), lower (upper) bound of the 95% confidence interval based on MCMC; $\text{L}_\text{N}$ ($\text{U}_\text{N}$), lower (upper) bound of the 95% confidence interval based on numerical derivatives.

boundaries for the two-sided 95% confidence interval for $\theta$ generated by numerical differentiation and MCMC. $L_M$ ($L_N$) and $U_M$ ($U_N$) denote the lower and upper bound of the confidence interval from the MCMC chain (numerical derivative). According to (17), Corollary 2 and Theorem 4, the terms $n|\text{MLE} - \text{CM}|$ ($n^{2/3}|\text{MLE} - \text{CM}|$), $\sqrt{n}|\text{SE}_M - \text{SE}_N|$ ($n^{1/6}|\text{SE}_M - \text{SE}_N|$), $n|L_M - L_N|$ ($n^{2/3}|L_M - L_N|$) and $n|U_M - U_N|$ ($n^{2/3}|U_M - U_N|$) for Cox regression with right censored data (current status data) in Table 1 (Table 2) are bounded in probability. And the realizations of these terms summarized in Tables 1 and 2 clearly illustrate their boundedness. Furthermore, we can conclude that the profile sampler based on the semiparametric models with faster convergence rate yields more accurate inferences about $\theta$.

5.3. *Partly linear normal model with current status data.* By differentiating the least favorable model with respect to $t$ or $\theta$, we can obtain

$$\dot{\ell}(t, \theta, k) = Q(x; t, k_t)(w - h_0(z)),$$

$$\ddot{\ell}(t, \theta, k) = (w - h_0(z))^2 \phi_t \left[ (1 - \delta) \frac{(1 - \Phi_t)q_t - \phi_t}{(1 - \Phi_t)^2} - \delta \frac{q_t \Phi_t + \phi_t}{\Phi_t^2} \right],$$

$$\ell_{t,\theta}(t, \theta, k) = (w - h_0(z))h_0(z)\phi_t \left[ (1 - \delta) \frac{(1 - \Phi_t)q_t - \phi_t}{(1 - \Phi_t)^2} - \delta \frac{q_t \Phi_t + \phi_t}{\Phi_t^2} \right],$$

$$\ell^{(3)}(t, \theta, k) = (w - h_0(z))^3 \phi_t R(q_t(x)),$$

$$\ell_{t,t,\theta}(t, \theta, k) = (w - h_0(z))^2 h_0(z) \phi_t R(q_t(x)),$$

$$\ell_{t,\theta,\theta}(t, \theta, k) = (w - h_0(z)) h_0^2(z) \phi_t R(q_t(x)),$$

where

$$R(q_t(x)) = \left[ (1 - \delta)\left( \frac{q_t^2 - 1}{1 - \Phi_t} + \frac{\phi_t q_t^2 - 2\phi_t q_t}{(1 - \Phi_t)^2} + \frac{2\phi_t^2}{(1 - \Phi_t)^3} \right) \right.$$
$$\left. - \delta \left( \frac{q_t^2 - 1}{\Phi_t} + \frac{3\phi_t q_t}{\Phi_t^2} + \frac{2\phi_t^2}{\Phi_t^3} \right) \right],$$

$q_t = q_{t,k_t(\theta,k)}(x)$, $\phi_t = \phi(q_t)$, and $\Phi_t = \Phi(q_t)$. The convergence rate for the estimated nuisance parameter is established in Lemma 4 by application of Theorem 1. The rate $r = 2/5$ is clearly faster than the cubic rate but slower than the parametric rate. Note that $\mathcal{O}_2^M$ is a $P$-Donsker class by technical tool T1 in the Appendix. Assumption C1 can be verified easily by recognizing that the three classes of functions specified in C1 depend on $(t, \theta, k)$ in a Lipschitz manner and are uniformly bounded. The remaining assumptions are verified in Lemmas 5 and 6 below.

LEMMA 4. *Under the above set-up for the partly linear normal model with current status data, we have for $\tilde{\theta}_n \xrightarrow{p} \theta_0$.*

$$(25) \qquad \|\hat{k}_{\tilde{\theta}_n} - k_0\|_{L_2} = O_P(n^{-2/5} + \|\tilde{\theta}_n - \theta_0\|).$$

LEMMA 5. *Under the above set-up for the partly linear normal model with current status data, assumptions* B1 *and* B2 *are satisfied.*

LEMMA 6. *Under the above set-up for the partly linear normal model with current status data, condition* (20) *is satisfied.*

**6. Future work.** It is clear that the estimation accuracy for $\theta$ in the profile sampler method is intrinsically determined by the semiparametric model specifications, specifically by the convergence rate of the nuisance parameter. Therefore it is very natural to raise a question about how to control the degree of accuracy. One potential strategy is to profile the penalized likelihood, whose penalty term is some norm on the nuisance parameter space such as the Sobolev norm. We expect that we can adjust the estimation accuracy of the proposed penalized profile sampler by tuning the corresponding smoothing parameter. We believe that under certain special model specifications, third or higher order semiparametric frequentist inference can be constructed by extending the Bartlett correction [1] and objective prior [24] results to semiparametric settings. There is a rich literature on the higher order properties of posteriors for parametric models and the choice of the prior; see, for example, [10, 19, 20].

## APPENDIX

PROOF OF THEOREM 2. We first prove (14). Note that

$$0 = \mathbb{P}_n \dot{\ell}(\hat{\theta}_n, \hat{\theta}_n, \hat{\eta}_n) = \mathbb{P}_n \dot{\ell}(\theta_0, \hat{\theta}_n, \hat{\eta}_n) + \mathbb{P}_n \ddot{\ell}(\theta_0, \hat{\theta}_n, \hat{\eta}_n)(\hat{\theta}_n - \theta_0)$$
$$+ \tfrac{1}{2}(\hat{\theta}_n - \theta_0)^T \otimes \mathbb{P}_n \ell^{(3)}(\theta_n^*, \hat{\theta}_n, \hat{\eta}_n) \otimes (\hat{\theta}_n - \theta_0),$$

where $\theta_n^*$ is in between $\theta_0$ and $\hat{\theta}_n$. By considering Lemma 2.1 below, we can derive the following:

$$0 = n^{-1} \sum_{i=1}^n \tilde{\ell}_0(x_i) + P\ddot{\ell}_0(\hat{\theta}_n - \theta_0) + n^{-1/2}\tilde{I}_0 \Delta_{4n}(\theta_0, \hat{\theta}_n, \hat{\eta}_n) \quad \text{and}$$

$$\sqrt{n}(\hat{\theta}_n - \theta_0) = \frac{1}{\sqrt{n}} \sum_{i=1}^n \tilde{I}_0^{-1} \tilde{\ell}_0(X_i) + \Delta_{4n}(\theta_0, \hat{\theta}_n, \hat{\eta}_n), \tag{26}$$

where

$$\Delta_{4n}(\theta_0, \hat{\theta}_n, \hat{\eta}_n) = \sqrt{n}\tilde{I}_0^{-1}\Delta_{1n}(\theta_0, \hat{\theta}_n, \hat{\eta}_n) + \sqrt{n}\tilde{I}_0^{-1}\Delta_{2n}(\theta_0, \hat{\theta}_n, \hat{\eta}_n)(\hat{\theta}_n - \theta_0)$$
$$+ \tfrac{1}{2}\sqrt{n}\tilde{I}_0^{-1}(\hat{\theta}_n - \theta_0)^T \otimes \mathbb{P}_n \ell^{(3)}(\theta_n^*, \hat{\theta}_n, \hat{\eta}_n) \otimes (\hat{\theta}_n - \theta_0)$$

and $\Delta_{1n}$ and $\Delta_{2n}$ are defined in the proofs of Lemma 2.1, respectively. The orders of magnitude of $\Delta_{1n}(\theta_0, \hat{\theta}_n, \hat{\eta}_n)$ and $\Delta_{2n}(\theta_0, \hat{\theta}_n, \hat{\eta}_n)$ obtained in the proofs of Lemma 2.1 imply that the order of magnitude of $\Delta_{4n}(\theta_0, \hat{\theta}_n, \hat{\eta}_n)$ is $O_P(M_n(r))$, as desired.

We next show (15). By (30) in Lemma 2.1 below, we have

$$
\begin{aligned}
\log pl_n(\hat{\theta}_n) &= \log pl_n(\theta_0) + (\hat{\theta}_n - \theta_0)^T \sum_{i=1}^n \tilde{\ell}_0(X_i) \\
&\quad - \frac{n}{2}(\hat{\theta}_n - \theta_0)^T \tilde{I}_0(\hat{\theta}_n - \theta_0) + \Delta_{3n}(\theta_0, \hat{\theta}_n, \hat{\eta}_n).
\end{aligned}
\tag{27}
$$

(30) also implies that

$$
\begin{aligned}
\log pl_n(\tilde{\theta}_n) &= \log pl_n(\hat{\theta}_n) + (\tilde{\theta}_n - \hat{\theta}_n)^T \left( \sum_{i=1}^n \tilde{\ell}_0(X_i) - n\tilde{I}_0(\hat{\theta}_n - \theta_0) \right) \\
&\quad - \frac{n}{2}(\tilde{\theta}_n - \hat{\theta}_n)^T \tilde{I}_0(\tilde{\theta}_n - \hat{\theta}_n) + \Delta_{3n}(\theta_0, \tilde{\theta}_n, \hat{\eta}_{\tilde{\theta}_n}) - \Delta_{3n}(\theta_0, \hat{\theta}_n, \hat{\eta}_n).
\end{aligned}
$$

Define $\Delta_{5n}(\tilde{\theta}_n, \hat{\theta}_n) = \log pl_n(\tilde{\theta}_n) - \log pl_n(\hat{\theta}_n) + (n/2)(\tilde{\theta}_n - \hat{\theta}_n)^T \tilde{I}_0(\tilde{\theta}_n - \hat{\theta}_n)$. By considering (26), we can obtain the respective upper and lower bounds of $\Delta_{5n}(\tilde{\theta}_n, \hat{\theta}_n)$ as follows:

$$
\begin{aligned}
\Delta_{5n}^U &= -\sqrt{n}(\tilde{\theta}_n - \hat{\theta}_n)^T \tilde{I}_0 \Delta_{4n}(\theta_0, \hat{\theta}_n, \hat{\eta}_n) \\
&\quad + \Delta_{3n}^U(\theta_0, \tilde{\theta}_n, \hat{\eta}_{\tilde{\theta}_n}) - \Delta_{3n}^L(\theta_0, \hat{\theta}_n, \hat{\eta}_n), \\
\Delta_{5n}^L &= -\sqrt{n}(\tilde{\theta}_n - \hat{\theta}_n)^T \tilde{I}_0 \Delta_{4n}(\theta_0, \hat{\theta}_n, \hat{\eta}_n) \\
&\quad + \Delta_{3n}^L(\theta_0, \tilde{\theta}_n, \hat{\eta}_{\tilde{\theta}_n}) - \Delta_{3n}^U(\theta_0, \hat{\theta}_n, \hat{\eta}_n),
\end{aligned}
$$

where $\Delta_{3n}^L$ and $\Delta_{3n}^U$ are defined in the proof of Lemma 2.1 and also shown to have magnitude $O_P(g_r(\|\tilde{\theta}_n - \hat{\theta}_n\|))$. Now the assumptions in Section 2 imply that $\Delta_{5n}^U(\tilde{\theta}_n, \hat{\theta}_n)$ and $\Delta_{5n}^U(\tilde{\theta}_n, \hat{\theta}_n)$ are of order $O_P(g_r(\|\tilde{\theta}_n - \hat{\theta}_n\|_2))$, and the proof is complete. $\square$

LEMMA 2.1. *Assuming the conditions of Theorem 2, we have*

$$
\mathbb{P}_n \dot{\ell}(\theta_0, \tilde{\theta}_n, \hat{\eta}_{\tilde{\theta}_n}) = \mathbb{P}_n \tilde{\ell}_0 + O_P(M_n(r) + \|\tilde{\theta}_n - \hat{\theta}_n\|)^2,
\tag{28}
$$

$$
\mathbb{P}_n \ddot{\ell}(\theta_0, \tilde{\theta}_n, \hat{\eta}_{\tilde{\theta}_n}) = P\ddot{\ell}_0 + O_P(M_n(r) + \|\tilde{\theta}_n - \hat{\theta}_n\|),
\tag{29}
$$

$$
\begin{aligned}
\log pl_n(\tilde{\theta}_n) &= \log pl_n(\theta_0) + (\tilde{\theta}_n - \theta_0)^T \sum_{i=1}^n \tilde{\ell}_0(X_i) \\
&\quad - \frac{n}{2}(\tilde{\theta}_n - \theta_0)^T \tilde{I}_0(\tilde{\theta}_n - \theta_0) + O_P(g_r(\|\tilde{\theta}_n - \hat{\theta}_n\|))
\end{aligned}
\tag{30}
$$

*for any random sequence $\tilde{\theta}_n - \hat{\theta}_n \xrightarrow{p} 0$.*

Proof.   By Taylor expansion of $\theta \mapsto P\dot{\ell}(\theta_0, \theta, \hat{\eta}_{\hat{\theta}_n})$, we obtain:

$$P\dot{\ell}(\theta_0, \tilde{\theta}_n, \hat{\eta}_{\hat{\theta}_n}) = P\dot{\ell}(\theta_0, \theta_0, \hat{\eta}_{\hat{\theta}_n}) + P\ell_{t,\theta}(\theta_0, \theta_0, \hat{\eta}_{\hat{\theta}_n})(\tilde{\theta}_n - \theta_0)$$
$$+ \tfrac{1}{2}(\tilde{\theta}_n - \theta_0)^T \otimes P\ell_{t,\theta,\theta}(\theta_0, \theta_1^*, \hat{\eta}_{\hat{\theta}_n}) \otimes (\tilde{\theta}_n - \theta_0)$$
$$\equiv \Delta_1(\theta_0, \tilde{\theta}_n, \hat{\eta}_{\hat{\theta}_n}),$$

where $\theta_1^*$ is intermediate between $\tilde{\theta}_n$ and $\theta_0$. The assumptions in Section [2] imply $\Delta_1(\theta_0, \tilde{\theta}_n, \hat{\eta}_{\hat{\theta}_n})$ has order $O_P(M_n(r) + \|\tilde{\theta}_n - \hat{\theta}_n\|)^2$. By writing $\mathbb{G}_n(\dot{\ell}(\theta_0, \tilde{\theta}_n, \hat{\eta}_{\hat{\theta}_n}) - \tilde{\ell}_0)$ as the summation of $\mathbb{G}_n(\dot{\ell}(\theta_0, \tilde{\theta}_n, \hat{\eta}_{\hat{\theta}_n}) - \dot{\ell}(\theta_0, \theta_0, \hat{\eta}_{\hat{\theta}_n}))$ and $\mathbb{G}_n(\dot{\ell}(\theta_0, \theta_0, \hat{\eta}_{\hat{\theta}_n}) - \tilde{\ell}_0)$, we obtain that the difference between $\mathbb{P}_n\dot{\ell}(\theta_0, \tilde{\theta}_n, \hat{\eta}_{\hat{\theta}_n})$ and $\mathbb{P}_n\tilde{\ell}_0$ is

$$\Delta_{1n}(\theta_0, \tilde{\theta}_n, \hat{\eta}_{\hat{\theta}_n}) = \Delta_1(\theta_0, \tilde{\theta}_n, \hat{\eta}_{\hat{\theta}_n}) + n^{-1/2}\mathbb{G}_n\ell_{t,\theta}(\theta_0, \theta_2^*, \hat{\eta}_{\hat{\theta}_n})(\tilde{\theta}_n - \theta_0)$$
$$+ n^{-1/2}\mathbb{G}_n(\dot{\ell}(\theta_0, \theta_0, \hat{\eta}_{\hat{\theta}_n}) - \dot{\ell}_0),$$

where $\theta_2^*$ is intermediate between $\tilde{\theta}_n$ and $\theta_0$. The order of magnitude of $\Delta_{1n}(\theta_0, \tilde{\theta}_n, \hat{\eta}_{\hat{\theta}_n})$ follows from the assumptions [6] and C1. This completes the proof of [28].

By similar analysis, we obtain

$$P\ddot{\ell}(\theta_0, \tilde{\theta}_n, \hat{\eta}_{\hat{\theta}_n}) = P\ddot{\ell}_0 + \Delta_2(\theta_0, \tilde{\theta}_n, \hat{\eta}_{\hat{\theta}_n})$$

and

$$\mathbb{P}_n\ddot{\ell}(\theta_0, \tilde{\theta}_n, \hat{\eta}_{\hat{\theta}_n}) = P\ddot{\ell}_0 + \Delta_{2n}(\theta_0, \tilde{\theta}_n, \hat{\eta}_{\hat{\theta}_n}),$$

where $\Delta_2(\theta_0, \tilde{\theta}_n, \hat{\eta}_{\hat{\theta}_n}) = (\tilde{\theta}_n - \theta_0)^T \otimes P\ell_{t,t,\theta}(\theta_0, \theta_3^*, \hat{\eta}_{\hat{\theta}_n}) + [P\ddot{\ell}(\theta_0, \theta_0, \hat{\eta}_{\hat{\theta}_n}) - P\ddot{\ell}_0]$ and $\Delta_{2n}(\theta_0, \tilde{\theta}_n, \hat{\eta}_{\hat{\theta}_n}) = \Delta_2(\theta_0, \tilde{\theta}_n, \hat{\eta}_{\hat{\theta}_n}) + n^{-1/2}\mathbb{G}_n\ddot{\ell}(\theta_0, \tilde{\theta}_n, \hat{\eta}_{\hat{\theta}_n})$, and where $\theta_3^*$ is in between $\theta_0$ and $\tilde{\theta}_n$. The assumptions in Section [2] now yield the desired order of magnitude in [29].

Next, we will show [30]. Note that

$$n^{-1}(\log pl_n(\tilde{\theta}) - \log pl_n(\theta_0)) = \mathbb{P}_n\ell(\tilde{\theta}_n, \tilde{\theta}_n, \hat{\eta}_{\hat{\theta}_n}) - \mathbb{P}_n\ell(\theta_0, \theta_0, \hat{\eta}_{\theta_0}).$$

The right-hand side of the above equation is bounded below and above by $\mathbb{P}_n(\ell(\tilde{\theta}_n, \tilde{\psi}_n) - \ell(\theta_0, \tilde{\psi}_n))$, where the lower and upper bounds separately correspond to $\tilde{\psi}_n = (\theta_0, \hat{\eta}_{\theta_0})$ and $(\tilde{\theta}_n, \hat{\eta}_{\tilde{\theta}_n})$. By applying a three-term Taylor expansion to both upper and lower bounds, we obtain the corresponding upper bound, $\Delta_{3n}^U(\theta_0, \tilde{\theta}_n, \hat{\eta}_{\hat{\theta}_n})$, and lower bound, $\Delta_{3n}^L(\theta_0, \tilde{\theta}_n, \hat{\eta}_{\hat{\theta}_n})$, for $\Delta_{3n}(\theta_0, \tilde{\theta}_n, \hat{\eta}_{\hat{\theta}_n})$ defined as follows:

$$\Delta_{3n}(\theta_0, \tilde{\theta}_n, \hat{\eta}_{\hat{\theta}_n}) \equiv \log pl_n(\tilde{\theta}) - \log pl_n(\theta_0) - (\tilde{\theta}_n - \theta_0)^T \sum_{i=1}^n \tilde{\ell}_0(X_i)$$
$$+ \tfrac{1}{2}n(\tilde{\theta}_n - \theta_0)^T\tilde{I}_0(\tilde{\theta}_n - \theta_0),$$

where

$$\Delta_{3n}^U(\theta_0, \tilde{\theta}_n, \hat{\eta}_{\tilde{\theta}_n})$$
$$= n(\tilde{\theta}_n - \theta_0)^T \Delta_{1n}(\theta_0, \tilde{\theta}_n, \hat{\eta}_{\tilde{\theta}_n})$$
$$\quad + \frac{n}{2}(\tilde{\theta}_n - \theta_0)^T \Delta_{2n}(\theta_0, \tilde{\theta}_n, \hat{\eta}_{\tilde{\theta}_n})(\tilde{\theta}_n - \theta_0)$$
$$\quad + \frac{n}{6}\sum_{i=1}^d \sum_{j=1}^d \sum_{k=1}^d \mathbb{P}_n \ell_{t_i, t_j, t_k}(\theta_4^*, \tilde{\theta}_n, \hat{\eta}_{\tilde{\theta}_n})(\tilde{\theta}_n - \theta_0)_i (\tilde{\theta}_n - \theta_0)_j (\tilde{\theta}_n - \theta_0)_k,$$

$$\Delta_{3n}^L(\theta_0, \tilde{\theta}_n, \hat{\eta}_{\tilde{\theta}_n})$$
$$= n(\tilde{\theta}_n - \theta_0)^T \Delta_{1n}(\theta_0, \theta_0, \hat{\eta}_{\theta_0})$$
$$\quad + \frac{n}{2}(\tilde{\theta}_n - \theta_0)^T \Delta_{2n}(\theta_0, \theta_0, \hat{\eta}_{\theta_0})(\tilde{\theta}_n - \theta_0)$$
$$\quad + \frac{n}{6}\sum_{i=1}^d \sum_{j=1}^d \sum_{k=1}^d \mathbb{P}_n \ell_{t_i, t_j, t_k}(\theta_5^*, \theta_0, \hat{\eta}_{\theta_0})(\tilde{\theta}_n - \theta_0)_i (\tilde{\theta}_n - \theta_0)_j (\tilde{\theta}_n - \theta_0)_k,$$

where $\theta_4^*$ and $\theta_5^*$ are in between $\theta_0$ and $\tilde{\theta}_n$. (28) and (29) yield the order of $\Delta_{3n}^U(\theta_0, \tilde{\theta}_n, \hat{\eta}_{\tilde{\theta}_n})$ and $\Delta_{3n}^L(\theta_0, \tilde{\theta}_n, \hat{\eta}_{\tilde{\theta}_n})$, which is $O_P(g_r(\|\tilde{\theta}_n - \hat{\theta}_n\|))$.  □

PROOF OF THEOREM 3.  Suppose that $F_n(\cdot)$ is the posterior profile distribution of $\sqrt{n}\varrho_n$ with respect to the prior $\rho(\theta)$, where the vector $\varrho_n$ is defined as $\tilde{I}_0^{1/2}(\theta - \hat{\theta}_n)$. Let the parameter set for $\varrho_n$ be $\Xi_n$. The whole proof of Theorem 3 can be briefly summarized in the following expression:

$$F_n(\xi) = \frac{\int_{\varrho_n \in (-\infty, n^{-1/2}\xi] \cap \Xi_n} \rho(\hat{\theta}_n + \tilde{I}_0^{-1/2}\varrho_n) \frac{pl_n(\hat{\theta}_n + \tilde{I}_0^{-1/2}\varrho_n)}{pl_n(\hat{\theta}_n)} \, d\varrho_n}{\int_{\varrho_n \in \Xi_n} \rho(\hat{\theta}_n + \tilde{I}_0^{-1/2}\varrho_n) \frac{pl_n(\hat{\theta}_n + \tilde{I}_0^{-1/2}\varrho_n)}{pl_n(\hat{\theta}_n)} \, d\varrho_n}.$$

Note that $d\varrho_n$ above is the short notation for $d\varrho_{n1} \times \cdots \times d\varrho_{nd}$. We first partition the parameter set $\Xi_n$ as $\{\Xi_n \cap \{\|\varrho_n\| > r_n\}\} \cup \{\Xi_n \cap \{\|\varrho_n\| \le r_n\}\}$. By choosing the proper order of $r_n$, we find the posterior mass in the first partition is of arbitrarily small order and the mass inside the second partition region can be approximated by a stochastic polynomial in powers of $n^{-1/2}$ with an error of order dependent on the convergence rate of the nuisance parameter. This general approach applies to both the denominator and numerator, yielding a quotient series that leads to the desired result.

Before giving the formal proof, we need two intermediate lemmas:

LEMMA 3.1. *Choose $r_n = o(n^{-1/3})$ such that $\sqrt{n} r_n \to \infty$. Under the conditions of Theorem 3, we have*

$$(31) \quad \int_{\|\varrho_n\| > r_n} \rho(\hat{\theta}_n + \tilde{I}_0^{-1/2} \varrho_n) \frac{pl_n(\hat{\theta}_n + \tilde{I}_0^{-1/2} \varrho_n)}{pl_n(\hat{\theta}_n)} d\varrho_n = O_P(n^{-M}),$$

*for any positive number $M$.*

PROOF. Fix $r > 0$. We have

$$\int_{\|\varrho_n\| > r} \rho(\hat{\theta}_n + \tilde{I}_0^{-1/2} \varrho_n) \frac{pl_n(\hat{\theta}_n + \tilde{I}_0^{-1/2} \varrho_n)}{pl_n(\hat{\theta}_n)} d\varrho_n$$

$$\leq I\{\Delta_n^r < -n^{-1/2}\} \exp(-\sqrt{n}) \int_\Theta \rho(\theta) \, d\theta + I\{\Delta_n^r \geq -n^{-1/2}\},$$

where $\Delta_n^r \equiv \sup_{\|\varrho_n\| > r} \Delta_n(\hat{\theta}_n + \tilde{I}^{-1/2} \varrho_n)$. By Lemma 3.2 in [4], $I\{\Delta_n^r \geq -n^{-1/2}\} = O_P(n^{-M})$ for any fixed $r > 0$. This implies that there exists a positive decreasing sequence $r_n = o(n^{-1/3})$ with $\sqrt{n} r_n \to \infty$ such that (31) holds. □

LEMMA 3.2. *Choose $r_n = o(n^{-1/3})$ such that $\sqrt{n} r_n \to \infty$. Under the conditions of Theorem 3, we have*

$$(32) \quad \begin{aligned} \int_{\|\varrho_n\| \leq r_n} &\left| \frac{pl_n(\hat{\theta}_n + \tilde{I}_0^{-1/2} \varrho_n)}{pl_n(\hat{\theta})} \rho(\hat{\theta}_n + \tilde{I}_0^{-1/2} \varrho_n) - \exp\left(-\frac{n}{2} \varrho_n^T \varrho_n\right) \rho(\hat{\theta}_n) \right| d\varrho_n \\ &= O_P(n^{-1/2} M_n(r)). \end{aligned}$$

PROOF. The posterior mass over the region $\|\varrho_n\| \leq r_n$ is bounded by

$$(*) \quad \int_{\|\varrho_n\| \leq r_n} \left| \frac{pl_n(\hat{\theta}_n + \tilde{I}_0^{-1/2} \varrho_n)}{pl_n(\hat{\theta}_n)} \rho(\hat{\theta}_n) - \exp\left(-\frac{n}{2} \varrho_n^T \varrho_n\right) \rho(\hat{\theta}_n) \right| d\varrho_n$$

$$(**) \quad \begin{aligned} + \int_{\|\varrho_n\| \leq r_n} &\left| \frac{pl_n(\hat{\theta}_n + \tilde{I}_0^{-1/2} \varrho_n)}{pl_n(\hat{\theta}_n)} \rho(\hat{\theta}_n + \tilde{I}_0^{-1/2} \varrho_n) \right. \\ &\left. - \frac{pl_n(\hat{\theta}_n + \tilde{I}_0^{-1/2} \varrho_n)}{pl_n(\hat{\theta}_n)} \rho(\hat{\theta}_n) \right| d\varrho_n. \end{aligned}$$

Using (15) when $r \geq 1/2$, we obtain

$$(*) = \int_{\|\varrho_n\| \leq r_n} \left[ \rho(\hat{\theta}_n) \exp\left(-\frac{n \varrho_n^T \varrho_n}{2}\right) \left| \exp(O_P(n \|\varrho_n\|^3 + n^{-1/2})) - 1 \right| \right] d\varrho_n$$

$$= n^{-1/2} \int_{\|u_n\| \leq \sqrt{n} r_n} \left[ \rho(\hat{\theta}_n) \exp\left( -\frac{u_n^T u_n}{2} \right) \right.$$
$$\left. \times \left| \exp(n^{-1/2}(\|u_n\|^3 + 1) O_P(1)) - 1 \right| \right] du_n$$

$$= n^{-1} \times O_P(1) \times \int_{\|u_n\| \leq \sqrt{n} r_n} \left[ \rho(\hat{\theta}_n) \exp\left( -\frac{u_n^T u_n}{2} \right) (\|u_n\|^3 + 1) \right] du_n$$

$$= O_P(n^{-1}),$$

where the second equality follows by replacing $\sqrt{n} \varrho_n$ with $u_n$, and the third equality follows from the fact that $|\exp(n^{-1/2}(\|u_n\|^3 + 1) O_P(1)) - 1| = O_P(1) n^{-1/2} \times (\|u_n\|^3 + 1)$, since $n^{-1/2} \|u_n\|^3 = o(1)$.

However, when $1/4 < r < 1/2$, we obtain

$$(*) = \int_{\|\varrho_n\| \leq r_n} \left[ \rho(\hat{\theta}_n) \exp\left( -\frac{n \varrho_n^T \varrho_n}{2} \right) \left| \exp(O_P(g_r(\|\varrho_n\|))) - 1 \right| \right] d\varrho_n$$

$$= \int_{\|\varrho_n\| \leq r_n} \left[ \rho(\hat{\theta}_n) \exp\left( -\frac{n \varrho_n^T \varrho_n}{2} \right) |\bar{g}_r(\|\varrho_n\|) \times O_P(1)| \right] d\varrho_n.$$

When $1/4 < r < 1/3$, $O_P(g_r(\|\varrho_n\|)) = O_P(n^{1-2r}\|\varrho_n\| + n^{-2r+1/2})$. Note that there exists a $\delta > 0$ such that $r_n = n^{2r-1-\delta}$ satisfying $r_n = o(n^{-1/3})$ with $\sqrt{n} r_n \to \infty$ for any $1/4 < r < 1/3$. Therefore $n^{1-2r}\|\varrho_n\| + n^{-2r+1/2} = o(1)$ when $r_n$ is taken equal to $n^{2r-1-\delta}$ for some $\delta > 0$. In this case, it implies that $\bar{g}_r(\|\varrho_n\|) = n^{1-2r}\|\varrho_n\| + n^{-2r+1/2}$ for $1/4 < r < 1/3$. However $\bar{g}_r(\|\varrho_n\|) = g_r(\|\varrho_n\|)$ for $1/3 \leq r < 1/2$ since $g_r(\|\varrho_n\|) = o(1)$ when $r$ is in this range. Combining with previous analyses, we have $(*) = O_P(n^{-2r})$ for $1/4 < r < 1/2$. Summarizing the above analysis, we have $(*) = O_P(n^{-1/2} M_n(r))$.

By the following analysis of $(**)$, we will be able to show that $(**) = O_P(n^{-1})$ for $r \geq 1/2$ since $\exp(O_P(n\|\varrho_n\|^3 + n^{-1/2})) = O_P(1)$ with $\|\varrho_n\| \leq r_n$:

$$(**) = \int_{\|\varrho_n\| \leq r_n} \left[ |\dot{\rho}(\theta_n^*)^T \tilde{I}_0^{-1/2} \varrho_n| \exp\left( -\frac{n}{2} \varrho_n^T \varrho_n + O_P(n\|\varrho_n\|^3 + n^{-1/2}) \right) \right] d\varrho_n$$

$$\leq M \int_{\|\varrho_n\| \leq r_n} \left[ \|\varrho_n\| \exp\left( -\frac{n}{2} \varrho_n^T \varrho_n \right) \right] d\varrho_n$$

$$\times \sup_{\|\varrho_n\| \leq r_n} \exp(O_P(n\|\varrho_n\|^3 + n^{-1/2})),$$

where $\theta_n^*$ is intermediate between $\hat{\theta}_n$ and $\hat{\theta}_n + \tilde{I}_0^{-1/2} \varrho_n$.

In the case that $1/4 < r < 1/2$, we have the same conclusion:

$$(**) = \int_{\|\varrho_n\| \leq r_n} \left[ |\dot{\rho}(\theta_n^*)^T \tilde{I}_0^{-1/2} \varrho_n| \exp\left( -\frac{n}{2} \varrho_n^T \varrho_n + O_P(g_r(\|\varrho_n\|)) \right) \right] d\varrho_n$$

$$\lesssim \int_{\|\varrho_n\| \le r_n} \left[ \|\varrho_n\| \exp\left(-\frac{n}{2}\varrho_n^T \varrho_n\right) \right] d\varrho_n$$

$$+ \int_{\|\varrho_n\| \le r_n} \|\varrho_n\| \exp\left(-\frac{n}{2}\varrho_n^T \varrho_n\right) |\exp(O_P(g_r(\|\varrho_n\|))) - 1| \, d\varrho_n$$

$$\le O_P(n^{-1}) + O_P(n^{-2r-1/2}) = O_P(n^{-1}),$$

where $\theta_n^*$ is an intermediate value between $\hat{\theta}_n$ and $\hat{\theta}_n + \tilde{I}_0^{-1/2}\varrho_n$. The last inequality follows from the analysis of $(*)$ when $1/4 < r < 1/2$. Hence we have proved that $(**) = O_P(n^{-1})$ for $r > 1/4$. This completes the proof. $\square$

Next we start the formal proof of Theorem 3. Note first that

$$\int_{\varrho_n \in \Xi_n} \left[ \rho(\hat{\theta}_n + \tilde{I}_0^{-1/2}\varrho_n) \frac{pl_n(\hat{\theta}_n + \tilde{I}_0^{-1/2}\varrho_n)}{pl_n(\hat{\theta}_n)} \right] d\varrho_n$$

$$= \int_{\{\|\varrho_n\| > r_n\} \cap \Xi_n} \left[ \rho(\hat{\theta}_n + \tilde{I}_0^{-1/2}\varrho_n) \frac{pl_n(\hat{\theta}_n + \tilde{I}_0^{-1/2}\varrho_n)}{pl_n(\hat{\theta}_n)} \right] d\varrho_n$$

$$+ \int_{\{\|\varrho_n\| \le r_n\} \cap \Xi_n} \left[ \rho(\hat{\theta}_n + \tilde{I}_0^{-1/2}\varrho_n) \frac{pl_n(\hat{\theta}_n + \tilde{I}_0^{-1/2}\varrho_n)}{pl_n(\hat{\theta}_n)} \right] d\varrho_n.$$

By Lemma 3.1, the first integral on the right is of order $O_P(n^{-1/2}M_n(r))$. The second integral on the right can be decomposed into the following summands:

$$\int_{\{\|\varrho_n\| \le r_n\} \cap \Xi_n} \left[ \rho(\hat{\theta}_n + \tilde{I}_0^{-1/2}\varrho_n) \frac{pl_n(\hat{\theta}_n + \tilde{I}_0^{-1/2}\varrho_n)}{pl_n(\hat{\theta})} - \exp\left(-\frac{n}{2}\varrho_n^T\varrho_n\right)\rho(\hat{\theta}_n) \right] d\varrho_n$$

$$+ \int_{\{\|\varrho_n\| \le r_n\} \cap \Xi_n} \left[ \exp\left(-\frac{n}{2}\varrho_n^T\varrho_n\right)\rho(\hat{\theta}_n) \right] d\varrho_n.$$

The first part in the above is bounded by $O_P(n^{-1/2}M_n(r))$ via Lemma 3.2. The second part equals

$$n^{-1/2}\rho(\hat{\theta}_n) \int_{\{\|u_n\| \le \sqrt{n}r_n\} \cap \sqrt{n}\Xi_n} e^{-u_n^T u_n/2} \, du_n$$

$$= n^{-1/2}\rho(\hat{\theta}_n) \int_{\mathbb{R}^d} e^{-u_n^T u_n/2} \, du_n + O(n^{-1/2}M_n(r)),$$

where $u_n = \sqrt{n}\varrho_n$. The above equality follows from the inequality $\int_x^\infty e^{-y^2/2} \, dy \le x^{-1}e^{-x^2/2}$ for any $x > 0$.

Consolidating the above analysis, we obtain

$$(33) \quad \int_{\varrho_n \in \Xi_n} \left[ \rho(\hat{\theta}_n + \tilde{I}_0^{-1/2} \varrho_n) \frac{pl_n(\hat{\theta}_n + \tilde{I}_0^{-1/2} \varrho_n)}{pl_n(\hat{\theta}_n)} \right] d\varrho_n$$

$$= n^{-1/2} \rho(\hat{\theta}_n)(2\pi)^{d/2} + O_P(n^{-1/2} M_n(r)),$$

and, by similar analysis, we also have

$$(34) \quad \int_{\varrho_n \in (-\infty, n^{-1/2}] \cap \Xi_n} \left[ \rho(\hat{\theta}_n + \tilde{I}_0^{-1/2} \varrho_n) \frac{pl_n(\hat{\theta}_n + \tilde{I}_0^{-1/2} \varrho_n)}{pl_n(\hat{\theta}_n)} \right] d\varrho_n$$

$$= n^{-1/2} \rho(\hat{\theta}_n) \times \int_{(-\infty, \xi_1] \times \cdots \times (-\infty, \xi_d]} e^{-y^T y/2} \, dy + O_P(n^{-1/2} M_n(r)).$$

The quotient of (33) and (34) generates the desired error rate for fixed $\xi$. Note, however, that the above conclusions are unchanged if $\xi$ is replaced by an arbitrary sequence $\{\xi_n\} \in \mathbb{R}^d$. Thus (21) follows, proving Theorem 3 in its entirety. $\square$

PROOF OF COROLLARY 1. From the proof of Theorem 3, we have

$$\tilde{P}_{\theta|\tilde{X}}(\sqrt{n} \tilde{I}_0^{1/2}(\theta - \hat{\theta}_n) \leq \xi)$$

$$= \frac{\int_{\varrho_n \in (-\infty, n^{-1/2}\xi] \cap \Xi_n} \rho(\hat{\theta}_n + \tilde{I}_0^{-1/2} \varrho_n) \frac{pl_n(\hat{\theta}_n + \tilde{I}_0^{-1/2} \varrho_n)}{pl_n(\hat{\theta}_n)} \, d\varrho_n}{\int_{\varrho_n \in \Xi_n} \rho(\hat{\theta}_n + \tilde{I}_0^{-1/2} \varrho_n) \frac{pl_n(\hat{\theta}_n + \tilde{I}_0^{-1/2} \varrho_n)}{pl_n(\hat{\theta}_n)} \, d\varrho_n}.$$

By differentiating both sides relative to $\xi$ and combining with (33), we obtain

$$f_n(\xi) = \frac{\rho(\hat{\theta}_n + \frac{\tilde{I}_0^{-1/2}\xi}{\sqrt{n}}) \frac{pl_n(\hat{\theta}_n + n^{-1/2}\tilde{I}_0^{-1/2}\xi)}{pl_n(\hat{\theta}_n)}}{(2\pi)^{d/2} \rho(\hat{\theta}_n) + O_P(M_n(r))}.$$

By analysis similar to the proof of Corollary 2 in [4], the numerator equals $\rho(\hat{\theta}_n) \exp(-\xi^T \xi/2) + O_P(M_n(r))$. This completes the proof. $\square$

PROOF OF COROLLARY 2. We only show (23) in what follows. Expression (24) can be shown similarly. The expansion in (23) is the quotient of two expansions of the form (33) and (34). We can show this as follows: First,

$$\tilde{E}_{\theta|x}(\varrho_n) = \frac{\int_{\varrho_n \in \Xi_n} \varrho_n \rho(\hat{\theta}_n + \tilde{I}_0^{-1/2} \varrho_n) \frac{pl_n(\hat{\theta}_n + \tilde{I}_0^{-1/2} \varrho_n)}{pl_n(\hat{\theta}_n)} \, d\varrho_n}{\int_{\varrho_n \in \Xi_n} \rho(\hat{\theta}_n + \tilde{I}_0^{-1/2} \varrho_n) \frac{pl_n(\hat{\theta}_n + \tilde{I}_0^{-1/2} \varrho_n)}{pl_n(\hat{\theta}_n)} \, d\varrho_n}.$$

The denominator is $n^{-1/2}(2\pi)^{d/2}\rho(\hat{\theta}_n) + O_P(n^{-1/2}M_n(r))$ by (33). Similarly, by the proof of Theorem 3 we know the numerator is a random vector of the order $O_P(n^{-2r-1/2} + n^{-3/2})$. This yields the desired conclusion. $\square$

PROOF OF THEOREM 4. By Lemma 4.1 in [4], we can easily show that $\kappa_{n\alpha} = \tilde{I}_0^{-1/2}z_\alpha + O_P(M_n(r))$ for any $\xi < \alpha < 1 - \xi$ and some choice of $z_\alpha$, where $\xi \in (0, 1/2)$. Note that $\kappa_{n\alpha}$ is not unique since the $\alpha$th quantile of a $d$-dimensional standard normal distribution, $z_\alpha$, is not unique when $d > 1$. The classical Edgeworth expansion implies that $P(n^{-1/2}\sum_{i=1}^n \tilde{I}_0^{-1/2}\tilde{\ell}_0(X_i) \leq z_\alpha + a_n(\alpha)) = \alpha$, where $a_n(\alpha) = O(n^{-1/2})$, for $\xi < \alpha < 1 - \xi$. This $a_n(\alpha)$ is thus uniquely determined for each fixed $z_\alpha$ since $\tilde{\ell}_0(X_i)$ has at least one absolutely continuous component. Let $\hat{\kappa}_{n\alpha} = \tilde{I}_0^{-1/2}z_\alpha + (\sqrt{n}(\hat{\theta}_n - \theta_0) - n^{-1/2}\sum_{i=1}^n \tilde{I}_0^{-1}\tilde{\ell}_0(X_i)) + \tilde{I}_0^{-1/2}a_n(\alpha)$. Then $P(\sqrt{n}(\hat{\theta}_n - \theta_0) \leq \hat{\kappa}_{n\alpha}) = \alpha$. Combining with (14), we obtain $\hat{\kappa}_{n\alpha} = \kappa_{n\alpha} + O_P(M_n(r))$. The uniqueness of $\hat{\kappa}_{n\alpha}$ follows from that of $a_n(\alpha)$ for each fixed $z_\alpha$, up to a term of order $O_P(M_n(r))$. $\square$

PROOF OF THEOREM 5. Under the assumptions of Theorem 5, we next show that $\chi_b^{n\alpha} = \chi_{d,\alpha}^2 + O_P(M_n(r))$ for $\xi < \alpha < 1 - \xi$, where $\xi \in (0, 1/2)$. It is sufficient to show that $\tilde{P}_{\theta|\tilde{X}}(PLR_b(\theta) \leq \chi_{d,\alpha}^2) = \alpha + O_P(M_n(r))$ by considering the form of $PLR_b(\theta)$ and (22). Based on the analysis in the proof of Theorem 2, the term $O_P(g_r(\|\tilde{\theta}_n - \hat{\theta}_n\|))$ in (15) is actually bounded above by $\Delta_{5n}^U(\tilde{\theta}_n, \hat{\theta}_n)$ and bounded below by $\Delta_{5n}^L(\tilde{\theta}_n, \hat{\theta}_n)$. Thus it yields the inequality that $n(\theta - \hat{\theta}_n)^T \tilde{I}_0(\theta - \hat{\theta}_n) - \Delta_{5n}^U(\theta, \hat{\theta}_n) \leq PLR_b(\theta) \leq n(\theta - \hat{\theta}_n)^T \tilde{I}_0(\theta - \hat{\theta}_n) - \Delta_{5n}^L(\theta, \hat{\theta}_n)$ such that we have constructed the upper bound and lower bound for $\tilde{P}_{\theta|\tilde{X}}(PLR_b(\theta) \leq \chi_{d,\alpha}^2)$.

We next show that the upper and lower bound matches asymptotically with $\alpha + O_P(M_n(r))$. Without loss of generality, we only consider its upper bound in what follows:

$$\tilde{P}_{\theta|\tilde{X}}(n(\theta - \hat{\theta}_n)^T \tilde{I}_0(\theta - \hat{\theta}_n) \leq \chi_{d,\alpha}^2 + \Delta_{5n}^U(\theta, \hat{\theta}_n))$$

$$\leq \tilde{P}_{\theta|\tilde{X}}(W_n) + \tilde{P}_{\theta|\tilde{X}}(\|\varrho_n\| > r_n)$$

$$\leq \tilde{P}_{\theta|\tilde{X}}(W_n) + \frac{\int_{\{\|\varrho_n\| > r_n\} \cap \Xi_n} \rho(\hat{\theta}_n + \tilde{I}_0^{-1/2}\varrho_n)\frac{pl_n(\hat{\theta}_n + \tilde{I}_0^{-1/2}\varrho_n)}{pl_n(\hat{\theta}_n)}\, d\varrho_n}{\int_{\Xi_n} \rho(\hat{\theta}_n + \tilde{I}_0^{-1/2}\varrho_n)\frac{pl_n(\hat{\theta}_n + \tilde{I}_0^{-1/2}\varrho_n)}{pl_n(\hat{\theta}_n)}\, d\varrho_n}$$

$$\leq \tilde{P}_{\theta|\tilde{X}}(W_n) + O_P(n^{-M}),$$

where $r_n = o(n^{1/3})$ with $\sqrt{n}r_n \to \infty$, $W_n = \{n\varrho_n^T\varrho_n \leq \chi_{d,\alpha}^2 + \Delta_{5n}^U(\theta, \hat{\theta}_n)\} \cap \{\|\varrho_n\| \leq r_n\}$, and $M$ is an arbitrary positive number. The third inequality above follows from Lemma 3.1 and (33).

We next study $\tilde{P}_{\theta|\tilde{X}}(W_n)$. Accordingly,

$$
\tilde{P}_{\theta|\tilde{X}}(W_n) = \frac{\int_{W_n} \rho(\hat{\theta}_n + \tilde{I}_0^{-1/2}\varrho_n)\frac{pl_n(\hat{\theta}_n + \tilde{I}_0^{-1/2}\varrho_n)}{pl_n(\hat{\theta}_n)}\,d\varrho_n}{\int_{\Xi_n} \rho(\hat{\theta}_n + \tilde{I}_0^{-1/2}\varrho_n)\frac{pl_n(\hat{\theta}_n + \tilde{I}_0^{-1/2}\varrho_n)}{pl_n(\hat{\theta}_n)}\,d\varrho_n}
$$

$$
= \frac{\int_{W_n}[\rho(\hat{\theta}_n + \tilde{I}_0^{-1/2}\varrho_n)\frac{pl_n(\hat{\theta}_n + \tilde{I}_0^{-1/2}\varrho_n)}{pl_n(\hat{\theta}_n)} - \rho(\hat{\theta}_n)\exp(-\frac{n}{2}\varrho_n^T\varrho_n)]\,d\varrho_n}{n^{-1/2}\rho(\hat{\theta}_n)(2\pi)^{d/2} + O_P(n^{-1/2}M_n(r))}
$$

$$
\quad + \frac{\int_{W_n} \rho(\hat{\theta}_n)\exp(-\frac{n}{2}\varrho_n^T\varrho_n)\,d\varrho_n}{n^{-1/2}\rho(\hat{\theta}_n)(2\pi)^{d/2} + O_P(n^{-1/2}M_n(r))}
$$

$$
= O_P(M_n(r)) + \frac{\int_{W_n} \rho(\hat{\theta}_n)\exp(-\frac{n}{2}\varrho_n^T\varrho_n)\,d\varrho_n}{n^{-1/2}\rho(\hat{\theta}_n)(2\pi)^{d/2} + O_P(n^{-1/2}M_n(r))}
$$

$$
= \frac{\int_{V_n} \rho(\hat{\theta}_n)\exp(-\frac{n}{2}\varrho_n^T\varrho_n)\,d\varrho_n + \int_{W_n - V_n} \rho(\hat{\theta}_n)\exp(-\frac{n}{2}\varrho_n^T\varrho_n)\,d\varrho_n}{n^{-1/2}\rho(\hat{\theta}_n)(2\pi)^{d/2} + O_P(n^{-1/2}M_n(r))}
$$

$$
\quad + O_P(M_n(r))
$$

$$
= \alpha + \frac{\int_{W_n - V_n} \rho(\hat{\theta}_n)\exp(-\frac{n}{2}\varrho_n^T\varrho_n)\,d\varrho_n}{n^{-1/2}\rho(\hat{\theta}_n)(2\pi)^{d/2} + O_P(n^{-1/2}M_n(r))} + O_P(M_n(r)),
$$

where $V_n = \{n\varrho_n^T\varrho_n \leq \chi^2_{d,\alpha}\}$. The third equality in the above follows from (33) and Lemma 3.2 in the proof of Theorem 3. We next study the fraction in the last equality above. It is easy to show that $\{W_n - V_n\} \subseteq \{W_n - V_n \cap \{\|\varrho_n\| \leq r_n\}\} \subseteq \{\chi^2_{d,\alpha} \leq n\varrho_n^T\varrho_n \leq \chi^2_{d,\alpha} + \Delta^U_{5n}(\theta, \hat{\theta}_n)\} \cap \{\|\varrho_n\| \leq r_n\} \equiv T_n$. By replacing $\sqrt{n}\varrho_n$ with $u_n$, $T_n$ can be reexpressed as $\{\chi^2_{d,\alpha} \leq u_n^T u_n \leq \chi^2_{d,\alpha} + \Delta^U_{5n}(\theta, \hat{\theta}_n)\} \cap \{\|u_n\| \leq \sqrt{n}r_n\}$.

We next consider the order of $\int_{W_n - V_n} \rho(\hat{\theta}_n)\exp(-n\varrho_n^T\varrho_n/2)\,d\varrho_n$ for $r$ in different ranges. For $r \geq 1/2$, $\Delta^U_{5n}(\theta, \hat{\theta}_n) = O_P(n^{-1/2} + n^{-1/2}\|u_n\|^3)$. Under the condition that $\|u_n\| \leq \sqrt{n}r_n$, $\Delta^U_{5n}(\theta, \hat{\theta}_n) = o_P(1)$. Hence any subsequence of $u_n$ contained in $T_n$ is not diverging. In this case, $\Delta^U_{5n}(\theta, \hat{\theta}_n) = O_P(n^{-1/2})$. In summary we have the following inequalities:

$$
\int_{W_n - V_n} \rho(\hat{\theta}_n)\exp\left(-\frac{n}{2}\varrho_n^T\varrho_n\right)d\varrho_n \leq \int_{T_n} \rho(\hat{\theta}_n)\exp\left(-\frac{n}{2}\varrho_n^T\varrho_n\right)d\varrho_n
$$

$$
\leq n^{-1/2}\int_{Q_n} \rho(\hat{\theta}_n)\exp\left(-\frac{u_n^T u_n}{2}\right)du_n
$$

$$
\leq O_P(n^{-1}),
$$

where $Q_n \equiv \{\chi^2_{d,\alpha} \leq u_n^T u_n \leq \chi^2_{d,\alpha} + O_P(n^{-1/2})\} \cap \{\|u_n\| \leq \sqrt{n} r_n\}$. Hence $\tilde{P}_{\theta|\bar{X}}(W_n) = \alpha + O_P(n^{-1/2})$ when $r \geq 1/2$.

Similar arguments will now be applied to the case $1/4 < r < 1/2$. It is sufficient to show that $\int_{W_n - V_n} \rho(\hat{\theta}_n) \exp(-\frac{n}{2} \varrho_n^T \varrho_n) d\varrho_n = O_P(n^{-2r})$ for $r$ in the above range. When $1/3 \leq r < 1/2$ and $\|u_n\| \leq \sqrt{n} r_n$, $\Delta_{5n}^U(\theta, \hat{\theta}_n)$ converges to zero in probability since $\Delta_{5n}^U(\theta, \hat{\theta}_n) = O_P(n^{-1/2}\|u_n\|^3 + n^{-r}\|u_n\|^2 + n^{1/2-2r}\|u_n\| + n^{-2r+1/2})$. Consequently, $\Delta_{5n}^U(\theta, \hat{\theta}_n)$ is $O_P(n^{1/2-2r})$ by the analysis we used for the case when $r \geq 1/2$. However, for $1/4 < r < 1/3$, $\Delta_{5n}^U(\theta, \hat{\theta}_n) = O_P(n^{-2r+1/2}\|u_n\| + n^{-2r+1/2})$. By making the same choice of $r_n$ used in the proof of Lemma 3.2, we have $\Delta_{5n}^U(\theta, \hat{\theta}_n) = o_P(1)$. Hence there does not exist a diverging subsequence of $u_n$ contained in $T_n$ when we choose this specific $r_n$. This implies that $\Delta_{5n}^U(\theta, \hat{\theta}_n) = O_P(n^{1/2-2r})$. In other words, $\Delta_{5n}^U(\theta, \hat{\theta}_n) = O_P(n^{1/2-2r})$ for $1/4 < r < 1/2$. This implies that $\int_{W_n - V_n} \rho(\hat{\theta}_n) \exp(-\frac{n}{2} \varrho_n^T \varrho_n) d\varrho_n = O_P(n^{-2r})$. The same arguments also apply to the lower bound of $\tilde{P}_{\theta|\bar{X}}(PLR_b(\theta) \leq \chi^2_{d,\alpha})$. Thus we have shown that $\chi^{n\alpha}_b = \chi^2_{d,\alpha} + O_P(M_n(r))$ for $\xi < \alpha < 1 - \xi$, where $\xi \in (0, 1/2)$.

If we can show that $\chi^{n\alpha}_f = \chi^2_{d,\alpha} + O_P(M_n(r))$, then the whole proof is complete. Combining (14) and (15) in Theorem 2, we can rewrite $PLR_f(\theta_0)$ as $n^{-1} \sum_{i=1}^n \tilde{\ell}_0(X_i)^T \tilde{I}_0^{-1} \sum_{i=1}^n \tilde{\ell}_0(X_i)^T + O_P(M_n(r))$. By classical Edgeworth expansion, we have $P(n^{-1/2} \sum_{i=1}^n \tilde{I}_0^{-1/2} \tilde{\ell}_0(X_i) \leq z_\alpha) = \alpha + O(n^{-1/2})$, which directly yields $P(n^{-1} \sum_{i=1}^n \tilde{\ell}_0(X_i)^T \tilde{I}_0^{-1} \sum_{i=1}^n \tilde{\ell}_0(X_i)^T \leq \chi^2_{d,\alpha} + O(n^{-1/2})) = \alpha$. Thus $\chi^{n\alpha}_f = \chi^2_{d,\alpha} + O_P(M_n(r))$. This completes the proof. $\square$

PROOF OF LEMMA 2. We first review some known results from [18] about the Cox model with current status data. For some constant $C$ and every $x$ (under the assumed regularity conditions), we have

$$|lik(\theta_0, \Lambda_0)(x) - lik(\theta_0, \Lambda)(x)| \leq C|\Lambda(y) - \Lambda_0(y)|,$$

(35)
$$|\dot{\ell}(\theta_0, \theta_0, \Lambda)(x) - \dot{\ell}(\theta_0, \theta_0, \Lambda_0)(x)| \leq C|\Lambda(y) - \Lambda_0(y)|,$$

$$|lik(\theta_0, \Lambda)(x) - lik(\theta_0, \Lambda_0)(x) - A_0(\Lambda - \Lambda_0)(x) lik(\theta_0, \Lambda_0)(x)|$$
$$\leq |\Lambda(y) - \Lambda_0(y)|^2,$$

where $A_0 = A_{\theta_0, \Lambda_0}$ and $A_{\theta, \Lambda}$ is the score operator for $\Lambda$ at $(\theta, \Lambda)$, for example, the Fréchet derivative of $\log p_{\theta, \Lambda}$ relative to $\Lambda$. Thus by the decomposition of $P\dot{\ell}(\theta_0, \theta_0, \Lambda)$ in what follows, we can show (9) holds with the $L_2$ norm on $\Lambda$:

$$P\dot{\ell}(\theta_0, \theta_0, \Lambda) = P\left[\frac{p_0 - p_{\theta_0, \Lambda}}{p_0}(\dot{\ell}(\theta_0, \theta_0, \Lambda) - \dot{\ell}_0)\right]$$
$$- P\dot{\ell}(\theta_0, \theta_0, \Lambda_0)\left[\frac{p_{\theta_0, \Lambda} - p_0}{p_0} - A_0(\Lambda - \Lambda_0)\right].$$

Next we can show (7) by the following inequality:

$$P(\ddot{\ell}(\theta_0, \theta_0, \Lambda) - \ddot{\ell}(\theta_0, \theta_0, \Lambda_0))$$
$$\leq P|N(\theta_0, \theta_0, \Lambda) - N(\theta_0, \theta_0, \Lambda_0)| + P|\dot{\ell}^2(\theta_0, \theta_0, \Lambda) - \dot{\ell}^2(\theta_0, \theta_0, \Lambda_0)|$$
$$\leq P|N(\theta_0, \theta_0, \Lambda) - N(\theta_0, \theta_0, \Lambda_0)| + C\|\Lambda - \Lambda_0\|_{L_2},$$

where $N(t, \theta, \Lambda) = (\partial^2 lik(t, \Lambda_t(\theta, \Lambda))/\partial t^2)/lik(t, \Lambda_t(\theta, \Lambda))$. The second inequality follows from the boundedness of $\dot{\ell}(t, \theta, \Lambda)$ and (35). We next proceed to derive an upper bound for $P|N(\theta_0, \theta_0, \Lambda) - N(\theta_0, \theta_0, \Lambda_0)|$:

$$P|N(\theta_0, \theta_0, \Lambda) - N(\theta_0, \theta_0, \Lambda_0)|$$
$$\lesssim P|\Lambda Q(x; \theta_0, \Lambda) - \Lambda Q(x; \theta_0, \Lambda_0)| + P|\Lambda - \Lambda_0| + P|w(\Lambda) - w(\Lambda_0)|$$
$$\quad + P|w(\Lambda)\Lambda - w(\Lambda_0)\Lambda_0|$$
$$\lesssim P|\Lambda Q(x; \theta_0, \Lambda) - \Lambda Q(x; \theta_0, \Lambda_0)| + \|\Lambda - \Lambda_0\|_{L_2}$$
$$\lesssim \|\Lambda - \Lambda_0\|_{L_2},$$

where $w(\Lambda) = \phi(\Lambda)h_{00} \circ \Lambda_0^{-1} \circ \Lambda$. Clearly, $w(\Lambda_0) = h_0$. Note that $w(\Lambda)$ can be expressed as $\Lambda \varsigma(\Lambda) v(\Lambda)$, where $\varsigma(\Lambda) = \phi(\Lambda)/\Lambda$ and $v(\Lambda) = h_{00} \circ \Lambda_0^{-1} \circ \Lambda$. Note that $\varsigma(\Lambda)$ and $v(\Lambda)$ are both assumed bounded and Lipschitz. Hence $P|w(\Lambda) - w(\Lambda_0)| \lesssim \|\Lambda - \Lambda_0\|_{L_2}$ and $P|\Lambda w(\Lambda) - \Lambda_0 w(\Lambda_0)| \lesssim \|\Lambda - \Lambda_0\|_{L_2}$. This explains the second inequality in the above. The inequality that $P|\Lambda Q(x; \theta_0, \Lambda) - \Lambda Q(x; \theta_0, \Lambda_0)| \lesssim \|\Lambda - \Lambda_0\|_{L_2}$ follows from the inequality that $|(u(e^u - e^v))/((e^u - 1)(e^v - 1))| \lesssim |u - v|$ given that $u \geq 0$ in some compact set and $v > 0$ in some compact set. Combining this with the previous analysis, we can conclude that (7) holds under the given assumptions for the Cox model with current status data. Similar techniques can be applied to the verification of (8). We omit the details.

Finally, we only need to check (6). Note that $\mathbb{G}_n(\dot{\ell}(\theta_0, \theta_0, \hat{\Lambda}_{\tilde{\theta}_n}) - \dot{\ell}_0)$ can be written as follows:

$$\mathbb{G}_n((\hat{\Lambda}_{\tilde{\theta}_n} - \Lambda_0)zQ(x; \theta_0, \Lambda_0)) - \mathbb{G}_n((w(\hat{\Lambda}_{\tilde{\theta}_n}) - w(\Lambda_0))Q(x; \theta_0, \Lambda_0))$$
$$(36) \qquad + \mathbb{G}_n(Q(x; \theta_0, \hat{\Lambda}_{\tilde{\theta}_n})(z\hat{\Lambda}_{\tilde{\theta}_n} - w(\hat{\Lambda}_{\tilde{\theta}_n})) - Q(x; \theta_0, \Lambda_0)(z\Lambda_0 - w(\Lambda_0)))$$
$$- \mathbb{G}_n(((z\hat{\Lambda}_{\tilde{\theta}_n} - w(\hat{\Lambda}_{\tilde{\theta}_n})) - (z\Lambda_0 - w(\Lambda_0)))Q(x; \theta_0, \Lambda_0)).$$

To verify (6), we need to make use of the following technical tools:

Clearly, by Lemma 5.13 in [22] we have $\mathbb{G}_n((\hat{\Lambda}_{\tilde{\theta}_n} - \Lambda_0)zQ(x; \theta_0, \Lambda_0)) = O_P(n^{-1/6} + \|\hat{\Lambda}_{\tilde{\theta}_n} - \Lambda_0\|_{L_2}^{1/2})$, since $\alpha = 1$ for monotone functions $\Lambda$. Then by the relation $O_P(\|\hat{\Lambda}_{\tilde{\theta}_n} - \Lambda_0\|_{L_2}^{1/2}) = O_P((|\tilde{\theta}_n - \theta_0| + n^{-1/3})^{1/2}) = O_P(n^{-1/6} + n^{1/6}|\tilde{\theta}_n - \theta_0|)$, we know that the first line of (36) satisfies (6). Note that the class of Lipschitz functions of $\Lambda$ also has the same upper bound (i.e., $\alpha = 1$) for the entropy with bracketing number by Theorem 2.7.11 in [23] and the inequality that $N(\varepsilon, \mathcal{G}, \|\cdot\|_{L_2}) \leq N_{[\cdot]}(2\varepsilon, \mathcal{G}, \|\cdot\|_{L_2})$. This now implies that $\mathbb{G}_n((w(\hat{\Lambda}_{\tilde{\theta}_n}) - w(\Lambda_0))Q(x; \theta_0, \Lambda_0)) = O_P(n^{-1/6} + n^{1/6}|\tilde{\theta}_n - \theta_0|)$ since $w(\Lambda)$ is Lipschitz in $\Lambda$. Similar arguments apply to the other lines in (36). □

PROOF OF LEMMA 3. The proof is analogous to that of Lemma 2 in [14]. □

PROOF OF LEMMA 4. We apply Theorem 1 with $m_{\theta, k} = \lambda(\theta, k)$, where $\lambda(\theta, k) \equiv \log lik(\theta, k)$, since $lik(\theta, k)$ is bounded away from zero and infinity for $(\theta, k) \in \Theta \times \mathcal{O}_2^M$. It suffices to show (10) provided both $P(\lambda(\theta, k_0) - \lambda(\theta_0, k_0)) \gtrsim -\|\theta - \theta_0\|^2$ and $P(\lambda(\theta, k) - \lambda(\theta_0, k_0)) \lesssim -d_\theta^2(k, k_0)$ hold. Note that the maximality of the point $(\theta_0, k_0)$ around the criterion function $(\theta, k) \mapsto P\lambda(\theta, k)$ implies that $P(\lambda(\theta, k_0) - \lambda(\theta_0, k_0)) \gtrsim -\|\theta - \theta_0\|^2$. By using the inequality $P \log(q/p) \leq -h^2(p, q)$ and the relationship between Kullback–Leibler divergence and squared Hellinger distance $h$, we can show that $P(\lambda(\theta, k) - \lambda(\theta_0, k_0)) \leq -\int (\sqrt{p_{\theta,k}} - \sqrt{p_0})^2 \, d\mu \leq -\|p_{\theta,k} - p_0\|_{L_2}^2$. Hence $d_\theta(k, k_0) = \|p_{\theta,k} - p_0\|_{L_2}$. Thus we only need to verify condition (11) by Lemma 1 to complete the whole proof.

Condition (13) in Lemma 1 trivially holds by considering the forms of $m_{\theta,k}$ and $d_\theta(k, k_0)$. By Theorem 1, we can show that $d_{\tilde{\theta}_n}(\hat{k}_{\tilde{\theta}_n}, k_0) = O_P(\delta_n + \|\tilde{\theta}_n - \theta_0\|)$ for any $\delta_n$ satisfying $K(\delta_n, \mathcal{S}_{\delta_n}, L_2(P)) \leq \sqrt{n}\delta_n^2$, where the function $K$ is defined in (12). In other words, we need to calculate the $\varepsilon$-bracketing entropy number for the class of functions $\mathcal{S}_{\delta_n}$. To achieve the desired rate (25), we only need to show $H_B(\varepsilon, \mathcal{S}_{\delta_n}, L_2(P)) \leq \varepsilon^{-1/2}$ based on the above discussions. Recall that $\mathcal{S}_{\delta_n} \equiv \{x \mapsto \lambda(\theta, k)(x) - \lambda(\theta_0, k_0)(x) : d_\theta(k, k_0) < \delta_n, \|\theta - \theta_0\| < \delta_n\}$. By considering Lemma 9.24 in [11], we only need to show that $H_B(\varepsilon, \mathcal{C}, L_2(P)) \leq \varepsilon^{-1/2}$, where $\mathcal{C} \equiv \{x \mapsto \lambda(\theta, k)(x) : J_2(k) + \|k - k_0\|_\infty \leq C_1, \|\theta - \theta_0\| \leq C_1\}$.

Now we consider $\mathcal{C}_1 \equiv \{q_{\theta,k}(x)/(1 + J(k)) : \|k - k_0\|_\infty \leq C_1, \|\theta - \theta_0\| \leq C_1\}$. By technical tool T1 below, we obtain $H_B(\varepsilon, \mathcal{C}_1, L_2(P)) \leq \varepsilon^{-1/2}$ as desired.

T1. (See [2].) For each $0 < C < \infty$ and $\delta > 0$ we have

$$(37) \qquad H_B(\delta, \{\eta : \|\eta\|_\infty \leq C, J_k(\eta) \leq C\}, \|\cdot\|_\infty) \lesssim \left(\frac{C}{\delta}\right)^{1/k}.$$

Continuing, note that $\lambda(\theta, k)(X)$ can be rewritten as:

$$(38) \qquad \Delta \log \Phi(q_{\tilde{\theta},k}A) + (1 - \Delta)\log(1 - \Phi(q_{\tilde{\theta},k}A)),$$

where $A = 1 + J(k)$ and $\bar{q}_{\theta,k} \in \mathcal{C}_1$. We next calculate the $\varepsilon$-bracketing entropy number with the $L_2$ norm for the class of functions $R_1 \equiv \{k_a(t) : t \mapsto \log \Phi(at)$ for $a \geq 1$ and $t \in T\}$, where $T$ is some bounded subset in $\mathbb{R}^1$. Note that $k_a(t)$ is increasing (decreasing) in $a$ for $t > 0$ $(t < 0)$. After some derivation, we obtain that $\sup_{t \in T} |k_a(t) - k_b(t)| \lesssim |a - b|$ for any fixed $a, b > 1$ and $\sup_{a, b \geq A_0, t \in T} |k_a(t) - k_b(t)| \lesssim A_0^{-1}$. The above two inequalities imply that the $\varepsilon$-bracketing number with uniform norm is of order $O(\varepsilon^{-2})$ for $a \in [1, \varepsilon^{-1}]$ and is 1 for $a > \varepsilon^{-1}$. Thus we know $H_B(\varepsilon, R_1, L_2) = O(\log \varepsilon^{-2})$. By applying a similar analysis to $R_2 \equiv \{k_a(t) : t \mapsto \log(1 - \Phi(at)) \text{ for } a \geq 1 \text{ and } t \in T\}$, we obtain that $H_B(\varepsilon, R_2, L_2) = O(\log \varepsilon^{-2})$. This combined with Lemma 15.2 in [11], yields that $H_B(\varepsilon, \mathcal{C}, L_2) \lesssim \varepsilon^{-1/2}$. Thus far we have shown that $d_{\tilde{\theta}_n}(\hat{k}_{\tilde{\theta}_n}, k_0) = O_P(n^{-2/5} + \|\tilde{\theta}_n - \theta_0\|)$. Now by the usual Taylor expansion and the assumption that $EVar(W|Z)$ is positive definite, we have verified (25). $\square$

PROOF OF LEMMA 5. Note that $k \in \mathcal{O}_2^M$. Hence we can easily verify assumption B1 since every map $(t, \theta, k) \mapsto (\partial^{l+m}/\partial t^l \partial \theta^m) \ell(t, \theta, k)$ is uniformly bounded. Note that $(C, W)$ lies in some bounded set and $h_0(\cdot)$ is bounded. Hence we can show that the Fréchet derivatives of $k \mapsto \dot{\ell}(\theta_0, \theta_0, k)$ and $k \mapsto \ell_{t, \theta}(\theta_0, \theta_0, k)$ for any $k \in \mathcal{O}_2^M$ are bounded operators, that is, $|\ddot{\ell}(\theta_0, \theta_0, k)(X) - \ddot{\ell}_0(X)|$ is bounded by the product of some integrable function and $|k - k_0|(Z)$. Thus (7) and (8) are satisfied, and the bounded Fréchet derivative of $k \mapsto \dot{\ell}(\theta_0, \theta_0, k)$ plus second-order Fréchet differentiability of $k \mapsto lik(\theta_0, k)$ implies (9).

Since the convergence rate $r = 2/5$, it suffices to show the asymptotic equicontinuity condition (6), provided (39) holds. Accordingly,

$$(39) \qquad \mathbb{G}_n(\dot{\ell}(\theta_0, \theta_0, \hat{k}_{\tilde{\theta}_n}) - \dot{\ell}_0) = O_P(n^{-3/10} + n^{1/10}\|\tilde{\theta}_n - \theta_0\|).$$

To show (39), we need the following technical tool T2:

T2. (Lemma 3.4.2 in [23].) Let $\mathcal{F}$ be a class of measurable functions such that $Pf^2 < \delta^2$ and $\|f\|_\infty \leq M$ for every $f$ in $\mathcal{F}$. Then

$$E_P^*\|\mathbb{G}_n\|_\mathcal{F} \lesssim K(\delta, \mathcal{F}, L_2(P))\left(1 + \frac{K(\delta, \mathcal{F}, L_2(P))}{\delta^2\sqrt{n}}M\right),$$

where $K(\delta, \mathcal{F}, \|\cdot\|) = \int_0^\delta \sqrt{1 + H_B(\varepsilon, \mathcal{F}, \|\cdot\|)}\,d\varepsilon$.

To utilize this tool, note first that (25) implies:

$$P\left(\frac{\dot{\ell}(\theta_0, \theta_0, \hat{k}_{\tilde{\theta}_n}) - \dot{\ell}_0}{n^{-3/10} + n^{1/10}\|\tilde{\theta}_n - \theta_0\|}\right)^2 \lesssim O_P(n^{-1/5}).$$

We next define the set $\mathcal{Q}_n$ as

$$\{g \in L_2(P) : Pg^2 \leq C_n n^{-1/5}\}$$

$$\cap\left\{\frac{\dot{\ell}(\theta_0, \theta_0, k) - \dot{\ell}_0}{n^{-3/10} + n^{1/10}\|\theta - \theta_0\|} : k \in \mathcal{O}_2^M, \|\theta - \theta_0\| \leq \delta\right\},$$

for some $\delta > 0$. Obviously the function $(\dot{\ell}(\theta_0, \theta_0, \hat{k}_{\tilde{\theta}_n}) - \dot{\ell}_0)/(n^{-3/10} + n^{1/10}\|\tilde{\theta}_n - \theta_0\|) \in \mathcal{Q}_n$ on a set of probability arbitrarily close to one, as $C_n \to \infty$. If we can show $\lim_{n\to\infty} E^*\|\mathbb{G}_n\|_{\mathcal{Q}_n} < \infty$ by T2, then Lemma 5 is proved.

Note that $\ell(\theta_0, \theta_0, k)$ depends on $k$ in a Lipschitz manner. Consequently we can bound $H_B(\varepsilon, \mathcal{Q}_n, L_2(P))$ by the product of some constant and $H(\varepsilon, \mathcal{R}_n, L_2(P))$, where $\mathcal{R}_n$ is defined as $\{G_n(k) : J(G_n(k)) \lesssim n^{3/10}, \|G_n(k)\|_\infty \lesssim n^{3/10}\}$, and where $G_n(k) = k/(n^{-3/10} + n^{1/10}\|\theta - \theta_0\|)$. By the main results in [2], we know $H(\varepsilon, \mathcal{R}_n, L_2(P)) \lesssim (n^{3/10}/\varepsilon)^{1/k}$. Note that $\delta_n = n^{-1/10}$ and $M_n = n^{3/10}$ in T2. Thus by calculation using T2, we establish $\lim_{n\to\infty} E^*\|\mathbb{G}_n\|_{\mathcal{Q}_n} < \infty$. $\square$

PROOF OF LEMMA 6. By the assumption that $\Delta_n(\tilde{\theta}_n) = o_P(1)$, we have $\Delta_n(\tilde{\theta}_n) - \Delta_n(\theta_0) \geq o_P(1)$. Thus the following inequality holds:

$$n^{-1}\sum_{i=1}^{n}\log\left[\frac{H(\tilde{\theta}_n, \hat{k}_{\tilde{\theta}_n}; X_i)}{H(\theta_0, \hat{k}_{\theta_0}; X_i)}\right] \geq o_P(1),$$

where $H(\theta, k; X) = \Delta\Phi(C - \theta W - k(Z)) + (1 - \Delta)(1 - \Phi(C - \theta W - k(Z)))$. By the assumptions on $k$, we know that $H(\tilde{\theta}_n, \hat{k}_{\tilde{\theta}_n}; X_i)$ belongs to some $P$-Donsker class. Combining the above conclusion and the inequality $\alpha\log x \leq \log(1 + \alpha\{x - 1\})$ for some $\alpha \in (0, 1)$ and any $x > 0$, we can show that

$$(40) \qquad P\log\left[1 + \alpha\left(\frac{H(\tilde{\theta}_n, \hat{k}_{\tilde{\theta}_n}; X_i)}{H(\theta_0, \hat{k}_{\theta_0}; X_i)} - 1\right)\right] \geq o_P(1).$$

The strict concavity of $x \mapsto \log(1 + \alpha(x - 1))$ ensures that

$$P\log\left[1 + \alpha\left(\frac{H(\tilde{\theta}_n, \hat{k}_{\tilde{\theta}_n}; X_i)}{H(\theta_0, \hat{k}_{\theta_0}; X_i)} - 1\right)\right] \leq 0.$$

This combined with (40) implies that

$$P\log\left[1 + \alpha\left(\frac{H(\tilde{\theta}_n, \hat{k}_{\tilde{\theta}_n}; X_i)}{H(\theta_0, \hat{k}_{\theta_0}; X_i)} - 1\right)\right] = o_P(1).$$

The strict concavity of $x \mapsto \log(1 + \alpha(x - 1))$ forces the result that $P|\Phi(C - \tilde{\theta}_n W - \hat{k}_{\tilde{\theta}_n}(Z)) - \Phi(C - \theta_0 W - \hat{k}_{\theta_0})(Z)| = o_P(1)$. The desired conclusion now follows from model identifiability. $\square$

Department of Statistical Science
Duke University
214 Old Chemistry Building
Durham, North Carolina 27708
USA
E-mail: chengg@stat.duke.edu

Department of Biostatistics
University of North Carolina at Chapel Hill
3101 McGavran–Greenberg Hall
Chapel Hill, North Carolina 27599
USA
E-mail: kosorok@unc.edu